\documentclass[12pt]{article}

\usepackage{tikz}
\usetikzlibrary{arrows.meta,bending,chains}
\usepackage{authblk}

\usepackage{lmodern}
\usepackage{lipsum}     
\usepackage{amsthm}
\usepackage{amssymb}
\usepackage{amsbsy}
\usepackage{stmaryrd}
\usepackage{mathtools}
\usepackage{bbm}
\usepackage{mathrsfs}
\usepackage{siunitx}    
\usepackage{pdflscape}  
\usepackage{rotating}   
\usepackage{textgreek}  
\usepackage{gensymb}    
\usepackage[misc]{ifsym} 
\usepackage{listings}   
\usepackage{colortbl}   
\usepackage{tabularx}   
\usepackage{longtable}  
\usepackage{subcaption}
\usepackage{multirow}
\usepackage{snotez}     
\usepackage{manfnt}
\usepackage{float}
\usepackage{booktabs,pifont}

\newcommand{\cmark}{\ding{51}}
\newcommand{\xmark}{\ding{55}}

\definecolor{myorange}{rgb}{0.9568,0.4941,0.1961}
\definecolor{myred}{rgb}{0.9098,0.1294,0.2078}
\definecolor{myblue}{rgb}{0.0352,0.4981,0.6509}
\definecolor{mygreen}{rgb}{0.2235,0.6353,0.2588}

\newcommand{\veps}{\varepsilon}
\newtheorem{remark}{Remark}

\newif\ifshowoptional
\showoptionalfalse 

\NewDocumentEnvironment{optional}{+b}
  {
    \ifshowoptional
      \par\begingroup\color{gray}#1\par\endgroup
    \fi
  }
  {}

\usepackage[commandnameprefix=always,todonotes={linecolor=myblue,bordercolor=myblue,backgroundcolor=white}]{changes}
\definechangesauthor[name=Phil, color=myred]{PM}

\title{Numerical homogenization of a linear multiscale conservation law}

\author{Claude Le Bris}
\author{Frédéric Legoll}
\author{Giulia Merlini}
\affil{\small École Nationale des Ponts et Chaussées, Institut Polytechnique de Paris, CNRS, 6 et 8 avenue Blaise Pascal, 77455 Marne-la-Vall\'ee, France}
\affil{\small MATHERIALS project-team, Inria Paris, 48 rue Barrault, 75013 Paris, France}
\affil{\small Emails: \{claude.le-bris,frederic.legoll,giulia.merlini\}@enpc.fr}

\begin{document}

\maketitle

\begin{abstract}
  We address the numerical approximation of a scalar conservation law with a highly oscillatory flux, in the linear setting. Taking advantage of the intrinsic numerical viscosity of the scheme employed that yields a parabolic regularization of the original hyperbolic equation, we use homogenization theory to derive an effective equation which in turn we solve numerically. Our numerical experiments show the performance of the approach on the example of the Lax-Friedrichs scheme.
\end{abstract}


\medskip



\section{Introduction}
This work addresses the numerical approximation of a multiscale scalar conservation law. In generality, the problem reads
\begin{equation}\label{eq: hyperbolic conservation law}
  \partial_t u^\veps(t,x) + \text{div} \, F \biggl( \frac{x}{\veps}, u^\veps(t,x) \biggr) = 0,
\end{equation}
and is complemented by suitable initial conditions and boundary conditions, along with the definition of the vector-valued, highly oscillatory flux $F$ depending upon the microscale $\veps$. All this will be made precise in the sequel. The numerical approximation of this problem is challenging because of the small-scale oscillations. In order to capture the expected oscillations of the scalar-valued solution~$u^\veps$, direct simulations indeed require resolving the microscopic scale $\veps$, thus leading to prohibitive computational costs.

\paragraph{Homogenization theory.} A classical approach to overcome this type of difficulty in the numerical approximation of solutions to partial differential equations consists in using homogenization theory (see \emph{e.g.}~\cite{bensoussan-lions-papanic,tartar2009intro,jikov,cioranescu-donato1999,blanc-bris2023} for a general exposition of the theory) to first identify the limit, as $\veps \to 0$, of the equation, the so-called \emph{homogenized equation}, which is in turn approximated numerically. Making the approach explicit and effective requires a \emph{geometric structure} of the oscillatory coefficients, typically the property of being periodic in the small scale. Variants of this strategy, known as numerical homogenization techniques, consist in using homogenization theory as only a guideline in order to design multiscale numerical methods. The case of oscillatory \emph{elliptic} equations (or systems) is the prototypical case where these strategies may be put in action. The case of \emph{parabolic} equations, close in nature to that of elliptic equations, is also a favorable ground (and we will actually return to this case in the sequel).

As will be recalled shortly, the \emph{hyperbolic} nature of the conservation law~\eqref{eq: hyperbolic conservation law} unfortunately creates a much more delicate context for homogenization theory and the numerical approaches that follow from this paradigm. In the present work, we focus on the case of a \emph{linear} flux
\begin{equation*}
  F \biggl( \frac{x}{\veps}, u^\veps(t,x) \biggr) = u^\veps(t,x) \, b\left( \frac{x}{\veps} \right),
\end{equation*}
with a velocity field $b$ that, for simplicity, we take smooth and periodic with respect to the small scale $\veps$. Then, \eqref{eq: hyperbolic conservation law} reduces to the linear transport equation in conservative form
\begin{equation}\label{eq: linear hyperbolic conservation law}
  \partial_t u^\veps(t,x) + \text{div} \, \biggl( u^\veps(t,x) \, b \biggl(\frac{x}{\veps}\biggr) \biggr) = 0.
\end{equation}
The linear setting~\eqref{eq: linear hyperbolic conservation law}, however simple it is, is already both practically relevant and sufficiently challenging to constitute a realistic test case for the approach we are going to develop. We will henceforth only consider this case in the present work.

In the one-dimensional case, it is classical to prove, using the method of lines and an analytic argument restricted to the one-dimensional setting, that the linear transport equation~\eqref{eq: linear hyperbolic conservation law} converges, as~$\veps \to 0$, to a linear transport equation of the form
\begin{equation} \label{eq:effective-1d-lin}
\partial_t u^\star(t,x) + \text{div} \, \big( u^\star(t,x) \, b^\star \big) = 0,
\end{equation}
with, as effective velocity, the harmonic average
\begin{equation}\label{eq: harmonic average b}
  b^\star = \langle b^{-1} \rangle^{-1},
\end{equation}
where, here and throughout the article, we denote by~$\langle f_{\rm per} \rangle$ the average of a periodic function~$f_{\rm per}$ over its periodic cell.

On the other hand, the problem of homogenization of~\eqref{eq: linear hyperbolic conservation law} becomes considerably more intricate in the two-dimensional setting. This setting essentially contains all the difficulties present in any arbitrary dimension. It is unclear that a homogenized equation exists, and there are indeed well-known cases where it does not. Challenging questions for the numerical approximation of~\eqref{eq: linear hyperbolic conservation law} arise from this theoretical observation. Somewhat more precisely, let us mention that, for general periodic velocity fields $b$, the equation indeed remains multiscale in the limit~$\veps \to 0$, and ``at best'' reads as a system of equations coupling the microscopic variable~$y$ (standing for~$x/\veps$ in the limit) and the macroscopic variable~$x$, which cannot be reduced to a single effective equation. Many questions in this topic remain open and are the subject of ongoing research, see~\cite{briane2020homogenization} for a recent work in this direction. In some specific situations however, the microscopic variable can be eliminated. This is the case, for instance, for \emph{ergodic flows}~\cite{hou1992homogenization, tassa1997homogenization}, for which the limit equation is a linear transport equation of the form~\eqref{eq:effective-1d-lin} only involving the macroscopic variable~$x$ and an effective velocity~$b^\star$. This is also the case for \emph{shear flows}~\cite{tartar1989,mascarenhas1993memory,weinan1992homogenization}, but then the limit macroscopic equation is a non-local diffusion equation with a memory term.

In sharp contrast to the purely hyperbolic case and in line with the elliptic case mentioned above, homogenization becomes clearer, in all dimensions, as soon as the original equation is slightly perturbed by a small viscous term adequately scaled to compete with the other terms in the limit~$\veps \to 0$ (see~\cite{dalibard2006homogenization}). This adequate scaling reads
\begin{equation}\label{eq: dalibard}
  \partial_t u^\veps(t,x) + \text{div} \, \biggl( u^\veps(t,x) \, b \biggl(\frac{x}{\veps}\biggr) \biggr) - \delta \, \veps \, \Delta_x u^\veps(t,x) = 0,
\end{equation}
where $\delta > 0$ is any fixed scalar parameter, say $\delta=1$ for simplicity. Irrespective of the velocity field and the geometric nature of the flow it generates, a homogenized limit then exists, in the form of a fully macroscopic equation, indeed hyperbolic in nature, since it is a linear transport equation of the type~\eqref{eq:effective-1d-lin}. Our reference study for this case is the work~\cite{dalibard2006homogenization} by A.-L. Dalibard, where the asymptotic analysis of the problem is carried out (including the case when the flux $F$ is quasilinear). We will often return to this reference in the present work since it serves as a theoretical guideline for our own work focused on computational approaches.

\paragraph{Our approach in short.} The present work is motivated by the previous observation. Rather than addressing the derivation of a possible homogenized limit of the continuous solution $u^\veps$ to~\eqref{eq: linear hyperbolic conservation law}, we consider a \emph{numerical scheme} approximating this equation (which is, in any event, all what one has access to, in most practically relevant cases). Since most numerical schemes contain a certain amount of viscosity, it is plausible (and we will establish below this is indeed the case) that, in an adequate regime of discretization parameters, we are able to derive a homogenized limit for the scheme itself (this notion will be made precise below). Consequently, this homogenized limit can in turn provide a practical strategy to approximate the discrete solution to the original scheme in a much more economical way than a direct utilization of the scheme.

In order to demonstrate the validity of our approach, we consider, to approximate~\eqref{eq: linear hyperbolic conservation law}, a scheme that is notoriously significantly diffusive, namely a Lax-Friedrichs (henceforth abbreviated as LxF) scheme. It reads, in the one-dimensional setting,
\begin{equation}\label{eq: lax-friedrich classic}
  \frac{U^{n+1}_i - \frac{1}{2} (U^n_{i+1} + U^n_{i-1})}{\delta t} + \frac{U^n_{i+1} \, b^\veps_{i+1} - U^n_{i-1} \, b^\veps_{i-1}}{2 h} = 0,
\end{equation}
with self-explanatory notations. The parameters $\delta t$ and $h$ define a uniform discretization in time and space, respectively. They are chosen according to the CFL condition
\begin{equation} \label{eq:cfl_cond_pre}
\delta t = \alpha \, h / \| b^\veps \|_{L^\infty} \quad \text{with $\alpha < 1$},
\end{equation}
to guarantee stability. From the scheme~\eqref{eq: lax-friedrich classic}, we will derive the \emph{equivalent} (a.k.a. \emph{modified}) equation
\begin{multline}\label{eq: equiv eq - lf intro}
  \partial_t u^\veps_{\rm mod} + \partial_x ( u^\veps_{\rm mod} \, b^\veps) + \frac{\delta t}{2} \biggl( \partial_x [ u^\veps_{\rm mod} \, b^\veps \, \partial_x b^\veps ] + (\partial_x [(b^\veps)^2]) \, \partial_x u^\veps_{\rm mod} \biggr) \\ - \left(\frac{h^2}{2\delta t} - \frac{\delta t \, (b^\veps)^2}{2} \right) \partial_{xx} u^\veps_{\rm mod} = 0,
\end{multline}
bearing in mind that all the details regarding the definition and the derivation of~\eqref{eq: equiv eq - lf intro}, which is quite classical in numerical analysis (see \emph{e.g.}~\cite{warming-hyett1974}), are somewhat delicate in our setting where the field $b^\veps = b(\cdot/\veps)$ oscillates. This will be clarified in the sequel. As expected for a scheme known to be viscous, the modified equation~\eqref{eq: equiv eq - lf intro} is parabolic in nature. Thanks to the CFL condition~\eqref{eq:cfl_cond_pre}, the prefactor of the second-order term indeed remains positive:
$$
\frac{h^2}{2\delta t} \left( 1 - \frac{\delta t^2 \, (b^\veps)^2}{h^2} \right) \geq \frac{h^2}{2\delta t} \, (1 - \alpha^2) > 0.
$$
For reasons that will be made clear below, we then adjust as follows the scaling of the two discretization parameters~$\delta t$ and~$h$ in function of the small parameter~$\veps$:
$$
h = c_x \, \veps, \qquad \delta t = c_t \, \veps, 
$$
with $\displaystyle c_t = \alpha \, c_x / \| b^\veps \|_{L^\infty}$ and $\alpha<1$ to satisfy the CFL condition. With this specific scaling of the discretization parameters, we thus obtain from~\eqref{eq: equiv eq - lf intro} a parabolic equation only parameterized by the small scale $\veps$, \emph{formally} resembling~\eqref{eq: dalibard} (upon freezing~$b^\veps$), for a certain particular calibration of the parameter~$\delta$, namely
\begin{equation*}
  \delta = \frac{c_x^2}{2 c_t} \, (1 - \alpha^2).
\end{equation*}
The general machinery established in~\cite{dalibard2006homogenization}, and which we adapt here, then allows us to identify the homogenized limit of this modified equation~\eqref{eq: equiv eq - lf intro}, which reads as
\begin{equation} \label{eq:bar}
  \partial_t \bar{u}(t,x) + \bar{b} \, \partial_x \bar{u}(t,x) = 0,
\end{equation}
for a constant effective velocity $\bar{b}$ defined through the solution to a periodic cell problem (see~\eqref{eq: cell problem w1} below, in the simplest, one-dimensional setting) parameterized by $c_x$ and $c_t$. Later in our work, we will denote $\bar{u}$ by $\bar{u}_{c_x}$ (and likewise $\bar{b}$ by $\bar{b}_{c_x}$) to emphasize the dependency upon $c_x$ ($c_t$ being itself a function of $c_x$ as fixed above). Practically speaking, it then remains to solve (this cell problem and) the latter equation numerically, by one's method of choice, to obtain an efficient approximation of the numerical solution that would have been provided by the original Lax-Friedrichs scheme put in action on the original equation~\eqref{eq: linear hyperbolic conservation law}, presumably using much finer discretization parameters. The point is, of course, that the small scale~$\veps$ having disappeared in~\eqref{eq:bar}, coarser discretization parameters in time and space can be adopted. We cannot emphasize enough this shift of focus: \emph{our object of reference is the original scheme, not the equation}. In layman's terms, we acknowledge that one has, for some reason, decided that the scheme to use is the Lax-Friedrichs scheme, but, given the oscillatory character of the velocity field, one cannot put it in action because it would require too small discretization parameters. Using the intrinsic viscous nature of the scheme and guided by homogenization theory, we suggest an alternative route that is likely to provide better performance. Our numerical experiments will show that it is indeed the case.

\medskip

Some comments are in order. 

\medskip

First, the approach is not restricted to the one-dimensional setting, which we chose above to summarize our strategy in the simplest possible manner. That said, the subtleties of theoretical homogenization for hyperbolic problems, present in higher dimensions and briefly recalled above, have their counterpart for our approach. They will be discussed at length in the sequel. Put simply, we will face the difficulty related to the fact that the two processes of \emph{homogenization} and \emph{discretization} do not necessarily commute, as is well known, and that this gets critical in the delicate framework of hyperbolic problems. There is indeed an intricate interplay between three ingredients: the discretization parameters $(\delta t,h)$, the prefactor~$\delta$ of the regularization term arising from the numerical viscosity of the scheme, and the microscale $\veps$ within the oscillatory flux. Depending on the relative scaling of these quantities, those we fix and those we vary, different asymptotic regimes arise, leading to significantly different behaviors of the numerical approximation.

Second, our choice of the Lax-Friedrichs scheme as starting point is motivated, as briefly mentioned above, by its well-known viscous character. It also turns out that, for this specific case, some technical details within our derivation are simpler. The Lax-Friedrichs scheme therefore provides a convenient framework for our analysis, allowing us to clearly identify and understand the effects of the interaction between discretization, regularization, and homogenization. Our approach could in principle be applied to any numerical scheme that contains viscosity, even if definite conclusions about the generality are yet to be obtained. In addition, our numerical experiments conclude that, practically, if we had replaced the Lax-Friedrichs scheme by another classical scheme, it would equally fail to be correct for a space discretization parameter~$h$ large as compared to $\veps$, and would perform equally well qualitatively when $h$ is small. From the practical perspective, the Lax-Friedrichs scheme is \emph{as good as any} for the specific purpose of demonstrating the superiority of the variant provided by our homogenization approach compared to the original version of the scheme.

Third, we have conducted our study in the setting of a \emph{linear} equation, with a \emph{periodic} oscillatory advection field. As already briefly mentioned above, the linear setting indeed provides a relevant and non-trivial starting point, and is already practically relevant \emph{per se} for many situations of the engineering and life sciences. Linear conservation equations are ubiquitous. The ultimate objective is, however, to address some nonlinear settings. As for our assumption regarding periodicity, which we confess may be seen as more artificial, it is of course motivated by our wish to proceed with homogenization theory while avoiding the theoretical technicalities (and the computational workload) of more general settings amenable to homogenization theory, such as the random setting. We are unsure the approach can be extended to those cases. So we have ongoing works on alternative strategies that will be reported on elsewhere.

A last remark we want to make in this introductory section is that the phenomena studied here have some similarities with the interesting questions investigated in~\cite{colombo2021} about the nonlocal-to-local limit for conservation laws in the context of traffic flow models. There also, discretization approaches interfere with the limit process at the continuous level, and some considerations analogous to ours are at play. Our perspective is however radically different, since our focus here is entirely computational and aims at constructing an actual numerical approximation of the equation considered.

\paragraph{Structure of our article.} Our article is organized as follows. The approach is described in details in Section~\ref{sec:approach}. To begin with, we derive in Section~\ref{subsec:modified} the modified equation for the Lax-Friedrichs scheme in the one-dimensional setting, and proceed with its homogenization in Section~\ref{subsec:homog} to obtain the effective equation. The derivation is extended to the two-dimensional setting in Section~\ref{subsec:extension2D}. The discretization of the effective equation is presented in Section~\ref{sec:generalities}, where we describe in particular the solution to the cell problem, so as to compute the effective velocity, along with some other important technical details such as the relevant regimes of discretization parameters and how they affect the results. The various metrics we will employ in our numerical tests to assess the quality of our results are also discussed.

Section~\ref{sec:numericresults} presents our numerical experiments. We specifically explore two asymptotic regimes:
\begin{itemize}
\item the regime where the discretization parameter~$c_x$ is kept fixed (and thus likewise $c_t$) while the microscale~$\veps$ of the numerical solution progressively vanishes (Section~\ref{subsec:numericseps1d}), thereby illustrating on the scheme itself the homogenization procedure of Section~\ref{sec:approach};
\item and the regime in which the microscale~$\veps$ is held fixed and small, while the discretization parameter~$c_x$ is taken finer and finer (Section~\ref{subsec:refinementregime}), so that the numerical solution of the original scheme approaches the exact, oscillatory solution of the original conservation law.
\end{itemize}
%
The latter regime is examined in dimension one (Section~\ref{subsubsec:refine1d}), and next in dimension two (Section~\ref{subsec:numerics2D}), where the topology of the flow comes into play. Our tests in particular illustrate the difference between the case when the homogenized equation of the original conservation law~\eqref{eq: linear hyperbolic conservation law} is a classical transport equation, and a prototypical case, namely the shear flow, when the situation is far more intricate. In the former case, our approach provides an efficient surrogate not only for the specific numerical scheme (here the Lax-Friedrichs scheme) originally considered, which is the whole point of our efforts, but also for any precise numerical solution of the hyperbolic conservation law as $\veps$ vanishes. In the latter case, it is observed that the numerical solution of the effective equation obtained in our approach is still an efficient alternative to the Lax-Friedrichs scheme provided we adequately adjust the discretization parameters, but behaves radically differently from the solution of the original conservation law in the limit~$\veps \to 0$. This is not unexpected given that, following homogenization theory, the regularized (parabolic) version and the original version of the conservation law have different homogenized limits. Put differently, the two limits $(c_t,c_x) \to (0,0)$ and $\veps \to 0$ do not commute with one another. We conclude the article, in Section~\ref{subsec:performance}, by assessing the performance of the approach in the one-dimensional setting: we compare it with a classical finite difference scheme applied to the original equation, both in terms of the accuracy reached at a given computational cost, and of the cost required to reach a given accuracy.


\medskip

We regard this study as groundwork for more general settings: locally periodic media, nonlinear fluxes, higher-order schemes, different scalings.

\section{Our approach} \label{sec:approach}

As briefly announced in the introduction, we start from a numerical scheme that approximates the conservation law. For the sake of generality and to simplify extensions to nonlinear cases in future works, we choose a finite difference scheme that can be rewritten in a conservative form, thus as a finite volume scheme. We refer to~\cite{leveque2002fv} for an overview on conservative schemes for hyperbolic problems. Using a Taylor expansion of the discrete quantities, we are going to derive the associated equivalent (a.k.a. modified) equation, thus a family, indexed by the various discretization and scale parameters, of partial differential equations approximating the numerical scheme. In the following, we retain the equivalent equation up to the first non-trivial correction, neglecting higher-order terms. For each specific scaling of the discretizing parameters~$(\delta t,h)$ in terms of the microscale~$\veps$, an equivalent equation parameterized by~$\veps$ is obtained. We then proceed with homogenization.

\medskip

For simplicity, we consider~\eqref{eq: linear hyperbolic conservation law} on a bounded domain (say $\Omega = (0,1)^d$), with periodic boundary conditions. We first introduce the approach in one dimension ($d=1$) before extending it to the two-dimensional setting. Let $I \in \mathbb{N} \setminus \{0 \}$, $h = 1/I$ such that the grid points $(x_i)_{0 \leq i \leq I}$ are defined by $x_i = i \, h$. We consider a time interval $[0,T]$, with the final time $T > 0$, and a partition $t_n = n \, \delta t$, $n = 0, \dots, N$, for $N = [T / \delta t] + 1$, where the notation $[\cdot]$ indicates the integer part. We also define $x_{i + 1/2} = x_i + h/2$.

The numerical approximation to~\eqref{eq: linear hyperbolic conservation law} computed with the \textit{Lax-Friedrichs} (LxF) scheme at point $(x_i, t_n)$ is denoted with $(U^\veps)_i^n$. In its finite difference form, the scheme applied to~\eqref{eq: linear hyperbolic conservation law} reads
\begin{equation}\label{eq: lax-friedrich}
  \frac{(U^\veps)^{n+1}_i - \frac{1}{2} ((U^\veps)^n_{i+1} + (U^\veps)^n_{i-1})}{\delta t} + \frac{(U^\veps)^n_{i+1} \, b^\veps_{i+1} - (U^\veps)^n_{i-1} \, b^\veps_{i-1}}{2 h} = 0,
\end{equation}
where $(U^\veps)^n_i \approx u^\veps(t_n,x_i)$ and $b^\veps_i = b(x_i / \veps)$. Notice that this scheme is a natural extension of the classical LxF scheme for a flux~$F(u^\veps)$ only depending on the unknown~$u^\veps$ to the case of a flux $F(x / \veps,u^\veps)$ highly oscillatory in the space variable, and taken linear in~$u^\veps$ here.

Notice also that the Lax-Friedrichs scheme can be written in the \textit{conservative} form
\begin{equation*}
  (U^\veps)^{n+1}_i = (U^\veps)^n_i - \frac{\delta t}{h} \, \biggl[ F_{i+1/2}^n - F_{i-1/2}^n \biggr],
\end{equation*}
where $F_{i \pm 1/2}^n$ is the numerical flux function, which can be seen as an approximation of the flux $F$ at $x = x_{i \pm 1/2}$ and at time $t_n$. For the Lax-Friedrichs scheme applied to the linear equation~\eqref{eq: linear hyperbolic conservation law}, the numerical flux reads
\begin{equation*}
  F_{i + 1/2}^n = \frac{h}{2\delta t} \, \biggl[ (U^\veps)^n_i - (U^\veps)^n_{i+1} \biggr] + \frac{1}{2} \biggl[ (U^\veps)^n_{i+1} \, b^\veps_{i+1} + (U^\veps)^n_i \, b^\veps_i \biggr].
\end{equation*}

The LxF scheme~\eqref{eq: lax-friedrich} is stable under the CFL condition
\begin{equation}\label{eq: cfl}
  \delta t < h / \| b^\veps \|_{L^\infty(\Omega)}.
\end{equation}
In the following, since this will be sufficient for our arguments and in order to fix the ideas, we take
\begin{equation}\label{eq:cfl-general}
\delta t = \alpha \, h/ \| b^\veps \|_{L^\infty(\Omega)} \quad \text{for some $\alpha < 1$}.
\end{equation}

\subsection{Modified equation of the Lax-Friedrichs scheme (1D setting)} \label{subsec:modified}

Our intent here is to derive, at~$\veps$ fixed, the modified equation of the Lax-Friedrichs scheme, in the one-dimensional setting first. We perform a Taylor expansion of the quantities involved in~\eqref{eq: lax-friedrich} with respect to both the time and the space variables. Let us assume that there exists a smooth function $u^\veps_{\rm mod}$ such that
$$
u^\veps_{\rm mod}(t_n, x_i) = (U^\veps)^n_i.
$$
Denoting by $u^\veps_{\rm mod}$ the value of the solution at $(t_n, x_i)$, we formally obtain
\begin{multline}\label{eq: dev 1}
 u^\veps_{\rm mod}(t_{n+1}, x_i) = u^\veps_{\rm mod} + (\delta t) \, \partial_t u^\veps_{\rm mod} \\ + \frac{(\delta t)^2}{2} \, \partial_{tt} u^\veps_{\rm mod} + \mathcal{O} \left((\delta t)^3 \, \partial_t^3 u^\veps_{\rm mod} \right),
\end{multline}
\begin{multline}\label{eq: dev 2}
u^\veps_{\rm mod}(t_n, x_{i \pm 1}) = u^\veps_{\rm mod} \pm h \, \partial_x u^\veps_{\rm mod} + \frac{h^2}{2} \, \partial_{xx} u^\veps_{\rm mod} \\ \pm \frac{h^3}{6} \, \partial_{xxx} u^\veps_{\rm mod}+ \mathcal{O} \left(h^4 \, \partial_x^4 u^\veps_{\rm mod} \right),
\end{multline}
and
\begin{multline} \label{eq: dev 3}
(u^\veps_{\rm mod} \, b^\veps) \big|_{(t^n, x_{i \pm 1})} = u^\veps_{\rm mod} \, b^\veps \pm h \, \partial_x (u^\veps_{\rm mod} \, b^\veps) \\ + \frac{h^2}{2} \, \partial_{xx} (u^\veps_{\rm mod} \, b^\veps) + \mathcal{O}\left(h^3 \, \partial_x^3 (u^\veps_{\rm mod} \, b^\veps) \right),
\end{multline}
%
with obvious notations.

We insert~\eqref{eq: dev 1}-\eqref{eq: dev 2}-\eqref{eq: dev 3} in the numerical scheme~\eqref{eq: lax-friedrich} and derive at second order
\begin{multline} \label{eq: taylor eq}
  \frac{1}{\delta t} \biggl[ \biggl( u^\veps_{\rm mod} + (\delta t) \, \partial_t u^\veps_{\rm mod} + \frac{(\delta t)^2}{2} \, \partial_{tt} u^\veps_{\rm mod} \biggr) - \biggl( u^\veps_{\rm mod} + \frac{h^2}{2} \, \partial_{xx} u^\veps_{\rm mod} \biggr) \biggr] \\ = - \partial_x (u^\veps_{\rm mod} \, b^\veps) + \mathcal{O}_\veps \bigg(h^2 + \frac{h^4}{\delta t} + (\delta t)^2 \bigg),
\end{multline}
where the remainder term $\mathcal{O}_\veps$ parametrically depends upon~$\veps$ in a way that, here and below, could be made precise (see Remark~\ref{rk:epsilon} on this issue).

At first order, the modified equation thus reads
\begin{equation*}
  \partial_t u^\veps_{\rm mod} + \partial_x (u^\veps_{\rm mod} \, b^\veps) = \mathcal{O}_\veps \bigg(\delta t + \frac{h^2}{\delta t} \bigg).
\end{equation*}
Taking the time derivative, we have
\begin{equation*}
  \partial_{tt} u^\veps_{\rm mod} = \partial_x \left[ b^\veps \, \partial_x (u^\veps_{\rm mod} \, b^\veps) \right] + \mathcal{O}_\veps \bigg(\delta t + \frac{h^2}{\delta t} \bigg).
\end{equation*}
We insert this in~\eqref{eq: taylor eq} and obtain, now at second order,
\begin{multline}\label{eq: modif sim eps}
  \partial_t u^\veps_{\rm mod} + \partial_x (u^\veps_{\rm mod} \, b^\veps) + \frac{\delta t}{2} \, \partial_x [b^\veps \, \partial_x(u^\veps_{\rm mod} \, b^\veps)] \\ - \frac{h^2}{2\delta t} \, \partial_{xx} u^\veps_{\rm mod} = \mathcal{O}_\veps \bigg(h^2 + \frac{h^4}{\delta t} + (\delta t)^2 \bigg).
\end{multline}

At this stage, we choose to work in the regime where the discretization parameters~$h$ and~$\delta t$ scale as
\begin{equation} \label{eq:cxct}
  h = c_x \, \veps \qquad \text{and} \qquad \delta t = c_t \, \veps
\end{equation}
in function of~$\veps$, for two fixed, presumably small scalar constants~$c_x$ and~$c_t$. This regime is in line with practice, since using the LxF scheme for an oscillatory equation such as~\eqref{eq: linear hyperbolic conservation law} requires to at least take~$h < \veps$. The scaling in $\delta t$ follows because of the CFL condition we want to impose. More precisely, we additionally choose, in the vein of~\eqref{eq:cfl-general},
\begin{equation} \label{eq: scale cfl}
  c_t = \alpha \, c_x / \| b^\veps \|_{L^\infty(\Omega)} \quad \text{for some fixed $\alpha<1$}.
\end{equation}
Throughout our work, both in the current theoretical developments and later on in our numerical tests, we are thus left with the only two free parameters~$c_x$ and $\veps$ and henceforth denote the function~$u^\veps_{\rm mod}$ by~$u^\veps_{c_x}$ to make explicit both dependencies. Inserting~\eqref{eq:cxct} in the left-hand side of~\eqref{eq: modif sim eps} and putting to zero the remainder on the right-hand side, we obtain the modified equation
\begin{equation}\label{eq: equiv eq}
  \partial_t u^\veps_{c_x} + \partial_x (u^\veps_{c_x} \, b^\veps) + \frac{c_t \, \veps}{2} \, \partial_x [b^\veps \, \partial_x (u^\veps_{c_x} \, b^\veps)] - \frac{c_x^2 \, \veps}{2c_t} \, \partial_{xx} u^\veps_{c_x} = 0.
\end{equation}
In order to check that this equation is indeed parabolic in nature, we recognize it as
\begin{equation} \label{eq: equiv eq-bis}
  \partial_t u^\veps_{c_x} + \partial_x \biggl( u^\veps_{c_x} \, B\left(\frac{x}{\veps}\right) \biggr) - \veps \, \partial_x \biggl( A\left(\frac{x}{\veps}\right) \, \partial_x u^\veps_{c_x}\biggr) = 0,
\end{equation}
where the coefficients $A$ and $B$ are \emph{periodic}, here scalar-valued, and defined as
\begin{equation} \label{eq:AB-1d}
A(y) = \frac{c_x^2}{2c_t} - \frac{c_t}{2} \, b^2(y),
\qquad
B(y) = b(y) + \frac{c_t}{2} \, b(y) \, \partial_y b(y).
\end{equation}
We henceforth assume both of them smooth for simplicity, although they could be much more general.
 
This modified equation~\eqref{eq: equiv eq-bis} reads as (the one-dimensional version of) the original oscillatory equation~\eqref{eq: linear hyperbolic conservation law}, with a modified transport field and an additional diffusive term, which explicitly depend on the discretization parameters and the microscale. As mentioned in the introduction, the CFL condition~\eqref{eq: cfl}, that is~\eqref{eq: scale cfl}, ensures that the diffusion~$A$ is positive and bounded away from zero, so that~\eqref{eq: equiv eq-bis} is parabolic, and is similar in nature to the equation~\eqref{eq: dalibard} addressed throughout the work~\cite{dalibard2006homogenization}. Note, however, that the diffusion in~\eqref{eq: equiv eq-bis} varies (and even oscillates) in space, but the approach of~\cite{dalibard2006homogenization}, performed for a constant diffusion, still applies \emph{mutatis mutandis}.

\begin{remark} \label{rk:epsilon}
In the derivation~\eqref{eq: taylor eq}-\eqref{eq: modif sim eps} above, we have not kept track of the dependency of the derivatives of~$u^\veps_{c_x}$ in terms of the small scale~$\veps$, since~$\veps$ is evidently fixed in our derivation of the modified equation. If we want to make the remainder term $\mathcal{O}_\veps$ explicit in this derivation, we may assume that, due to the oscillatory nature of the velocity field $b^\veps(x) = b(x/\veps)$, the solution~$u^\veps_{c_x}$ itself oscillates at the microscale $\veps$, with, at most,
\begin{equation*} 
  \partial_t^k u^\veps_{c_x} = \mathcal{O}(\veps^{-k}) \qquad \text{and} \qquad \partial_x^k u^\veps_{c_x} = \mathcal{O}(\veps^{-k}).
\end{equation*}
This assumption is consistent with the argument of homogenization theory that will be performed in the next section, and turns out to be correct \emph{a posteriori}. In that case, simple calculations we skip here show that, instead of~$\displaystyle \mathcal{O}_\veps \bigg( h^2 + \frac{h^4}{\delta t} + (\delta t)^2 \bigg)$, the remainder on the right-hand side of~\eqref{eq: modif sim eps} actually reads~$\displaystyle \mathcal{O} \bigg( \frac{h^2}{\veps^3} + \frac{h^4 }{\delta t \, \veps^4} + \frac{(\delta t)^2}{\veps^3} \bigg)$, that is, given~\eqref{eq:cxct},
\begin{equation*} 
  \mathcal{O} \bigg( \frac{c_x^2}{\veps} + \frac{c_x^4 }{ c_t \, \veps} + \frac{c_t^2}{\veps} \bigg).
\end{equation*}
Consequently, the modified equation~\eqref{eq: equiv eq} indeed accurately approximates the numerical solution of the LxF scheme when, at~$\veps$ fixed, the discretization parameters~$c_x$ and~$c_t$ both vanish. However, the ``distance'' between this modified equation and the original LxF scheme might not go to zero as~$\veps$ itself vanishes. We will return to this when we comment upon our various numerical experiments in different settings and regimes of parameters. 
\end{remark}

\subsection{Homogenization of the modified equation} \label{subsec:homog}

We now perform the homogenization limit, as~$\veps \to 0$, of the modified equation~\eqref{eq: equiv eq} of the Lax-Friedrichs scheme. We assume that $u^\veps_{c_x}$ can be written as the two-scale expansion
\begin{equation}\label{eq: two-scale Ansatz}
  u^\veps_{c_x} (t,x) \approx u^0_{c_x} \left(t,x,\frac{x}{\veps}\right) + \veps \, u^1_{c_x} \left(t,x,\frac{x}{\veps}\right) + \dots,
\end{equation}
where each term $u^i_{c_x}(t,x,y)$ is periodic in the microscopic variable~$y$ within the periodic cell~$Y = (0,1)$. By substituting the Ansatz~\eqref{eq: two-scale Ansatz} into~\eqref{eq: equiv eq} and identifying the powers of $\veps$, we obtain a cascade of equations, the first two of which respectively read
\begin{align}
  \text{order $\veps^{-1}$:} & \quad \partial_y (b \, u^0_{c_x}) + \frac{c_t}{2} \, \partial_y \left[ b \, \partial_y (b \, u^0_{c_x}) \right] - \frac{c_x^2}{2c_t} \, \partial_{yy} u^0_{c_x} = 0,
  \label{eq: order -1}
  \\[6pt]
  \text{order $\veps^0$:} & \quad \partial_t u^0_{c_x} + \partial_x (b \, u^0_{c_x}) + \partial_y (b \, u^1_{c_x}) + \frac{c_t}{2} \left[ \partial_x (b \, \partial_y (b \, u^0_{c_x}) ) + \partial_y (b \, \partial_x (b \, u^0_{c_x}) ) \right]
  \nonumber
  \\
  & \hspace{1cm} + \frac{c_t}{2} \, \partial_y \left[ b \, \partial_y (b \, u^1_{c_x}) \right] - \frac{c_x^2}{c_t} \, \partial_{xy} u^0_{c_x} - \frac{c_x^2}{2c_t} \, \partial_{yy} u^1_{c_x} = 0.
  \label{eq: order 0} 
\end{align}
We next look for the dominant term $u^0_{c_x}$ under the general form
$$
u^0_{c_x}(t,x,y) = v(y, \bar{u}(x,t)),
$$ 
with $\displaystyle \langle u^0_{c_x}(t,x,\cdot) \rangle_Y = \bar{u}(t,x)$ and where the function $y \in Y \mapsto v(y,p)$ belongs to $H^1_{\rm per}(Y)$ for any $p \in \mathbb{R}$ (we recall that $\langle \cdot \rangle_Y$ is the average over the periodic cell $Y$). Given the linearity of our context, we expect $v$ to be linear with respect to its second variable and thus readily postulate that~$v(y,p) = p \, w(y)$, for some function~$w$ yet to be determined. We obtain from~\eqref{eq: order -1} that the \emph{corrector function}~$w$ must be solution to the \textit{cell problem}
\begin{equation}\label{eq: cell problem w1}
  \begin{cases}
    \partial_y (w \, B) - \partial_y ( A \, \partial_y w) = 0 \quad \text{in $\mathbb{R}$},
    \\[8pt]
    \text{$w$ is $Y$-periodic}, \qquad \langle w \rangle_Y = 1,
  \end{cases}
\end{equation}
with $A$ and $B$ defined in~\eqref{eq:AB-1d}. We note in passing that this problem is elliptic in nature, thus well posed, for the exact same reason as~\eqref{eq: equiv eq-bis} is parabolic. This is the key difference with the cell problem obtained for a purely hyperbolic equation, that reads as a first-order boundary value problem (possibly ill-posed in dimensions higher than one; we will return to this in Section~\ref{subsec:extension2D} below). In any event, it is evident to see, in this one-dimensional setting, that~\eqref{eq: cell problem w1} is well posed: it is even analytically solvable.

Averaging over~$y$ in~\eqref{eq: order 0} yields the \textit{LxF-homogenized problem}: $\bar{u}$ is the unique solution to the linear transport equation
\begin{equation}\label{eq: homog. problem LF}
  \begin{cases}
    \partial_t \bar{u}(t,x) + \bar{b} \, \partial_x \bar{u}(t,x) = 0,
    \\[4pt]
    \bar{u}(t = 0,x) = \bar{u}_0(x),
  \end{cases} 
\end{equation}
where the \textit{effective velocity} $\bar{b}$ is defined as
\begin{equation}\label{eq: effective velocity}
  \bar{b} = \int_Y w \, b + \frac{c_t}{2} \, b \, \partial_y (w \, b).
\end{equation}
Note that we have presented~\eqref{eq: homog. problem LF} and~\eqref{eq: effective velocity} using compact notations for $\bar{u}$ and $\bar{b}$. When needed, we will denote them by $\bar{u}_{c_x}$ and $\bar{b}_{c_x}$, to emphasize the dependency on $c_x$ ($c_t$ being itself a function of $c_x$ as given by~\eqref{eq: scale cfl}).

\medskip

In order for these formal manipulations to be correct and indeed give the limit of~\eqref{eq: equiv eq} as~$\veps \to 0$, we need to also adjust the initial condition we provide the equation with. As in~\cite{dalibard2006homogenization} and as in all homogenization results for hyperbolic equations themselves, the oscillatory equation~\eqref{eq: equiv eq} should be complemented by the \textit{well-prepared} initial condition
\begin{equation}\label{eq: well-prep}
  u^\veps_{c_x}(t = 0, x) = \bar{u}_0(x) \, w_{c_x}\left(\frac{x}{\veps}\right),
\end{equation}
where $\bar{u}_0$ is the initial condition in~\eqref{eq: homog. problem LF}. In~\eqref{eq: well-prep}, the oscillations in space at the small-scale reproduce those of the corrector function~$w$ solution to~\eqref{eq: cell problem w1}, a property that is key to establish strong convergence of oscillatory solutions in the limit $\veps \to 0$. For clarity, we have explicitly mentioned in~\eqref{eq: well-prep} the dependency of the corrector $w$ upon $c_x$.

\begin{remark}
In the absence of~\eqref{eq: well-prep}, it should be possible, as is done in~\cite[Theorem~1]{dalibard2007initial} for the continuous equation itself, to consider for our developments the following much less restrictive assumption: instead of~\eqref{eq: well-prep}, the initial condition reads $\displaystyle u^\veps_{c_x}(t = 0, x) = u_0\left(x,\frac{x}{\veps}\right)$ for a function $u_0(x,y)$ such that there exists two scalars~$\alpha$ and~$\beta$ such that, almost everywhere in~$(x,y) \in \Omega \times Y$,
\begin{equation*} 
  \alpha \, w(y) \leq u_0(x,y) \leq \beta \, w(y),
\end{equation*}
in which case, given the (small) parabolic nature of the problem studied, a boundary layer exists for small times. We do not proceed in this direction and prefer to concentrate ourselves on the case~\eqref{eq: well-prep} of well-prepared initial conditions.
\end{remark}

\subsection{Extension to the two-dimensional setting} \label{subsec:extension2D}

In the two-dimensional framework, we consider the oscillatory equation~\eqref{eq: linear hyperbolic conservation law} over a square $\Omega = (0,1)^2$, with periodic boundary conditions. To approximate $u^\veps$, we introduce a uniform grid $(x_{i,j})_{0 \leq i,j \leq I}$ where $x_{i,j} = (i \, h, j \, h)$. The velocity field $\displaystyle b^\veps = \left( p^\veps, q^\veps \right)$ evaluated at grid points is denoted with components $p^\veps_{i,j} = p(x_{i,j}/\veps)$ and $q^\veps_{i,j} = q(x_{i,j}/\veps)$. We work in the same regime~\eqref{eq:cxct} of discretization parameters and our argument proceeds exactly as in the one-dimensional setting. We start from the Lax-Friedrichs scheme
\begin{multline}\label{eq: lax-friedrich 2d}
  \frac{1}{\delta t} \biggl( (U^\veps)^{n+1}_{i,j} - \frac{1}{4} \biggl[ (U^\veps)^n_{i+1, j} + (U^\veps)^n_{i-1, j} + (U^\veps)^n_{i, j+1} + (U^\veps)^n_{i, j-1} \biggr] \biggr) \\ + \frac{1}{2h} \biggl( p^\veps_{i+1,j} \, (U^\veps)^n_{i+1,j} - p^\veps_{i-1,j} \, (U^\veps)^n_{i-1,j} \biggr) \\ + \frac{1}{2h} \biggl( q^\veps_{i,j+1} \, (U^\veps)^n_{i,j+1} - q^\veps_{i,j-1} \, (U^\veps)^n_{i,j-1} \biggr) = 0.
\end{multline}
We then derive the modified equation
\begin{equation} \label{eq: equiv eq-bis-2d}
  \partial_t u^\veps_{c_x} + \text{div} \left(u^\veps_{c_x} \, B\left(\frac{x}{\veps}\right)\right) - \veps \, \text{div} \left( A\left(\frac{x}{\veps}\right) \nabla u^\veps_{c_x} \right) = 0,
\end{equation}
with
\begin{equation} \label{eq:AB-2d}
A = \frac{c_x^2}{4 \, c_t} \, {\rm Id} - \frac{c_t}{2} \, b \otimes b,
\qquad
B = b + \frac{c_t}{2} \, (\text{div} \, b) \, b,
\end{equation}
(where ${\rm Id}$ denotes the identity matrix), which is the generalization of the one-dimensional formulae~\eqref{eq: equiv eq-bis}-\eqref{eq:AB-1d} to higher dimensional settings. The matrix $A$ is positive definite because of the CFL condition~\eqref{eq: scale cfl}, and bounded because of the regularity assumed on~$b$.

\medskip

Proceeding with the two-scale expansion and adapting our manipulations of the one-dimensional setting, we obtain likewise the two-dimensional cell problem
\begin{equation}\label{eq: cell problem 2d}
  \begin{cases}
    \text{div} \, (w \, B) - \text{div} \, (A \, \nabla w) = 0 \quad \text{in $\mathbb{R}^2$},
    \\[5pt]
    \text{$w$ is $Y$-periodic}, \qquad \langle w \rangle_Y = 1,
  \end{cases}
\end{equation}
with $A$ and $B$ of course given by~\eqref{eq:AB-2d}. The ellipticity of the cell problem is once again guaranteed by the CFL condition~\eqref{eq: scale cfl}. The well-posedness of~\eqref{eq: cell problem 2d} then follows from a classical argument based on Fredholm's alternative, as given in~\cite[Theorem~1.1]{droniou2009}. This is in sharp contrast with the case when the diffusion term is absent, leaving then the cell problem as a first order problem of the form~$\text{div} \, (w \, b) = 0$ which, depending upon the geometric properties of the flow associated to~$b$, might be, or not, well posed. The numerical viscosity embedded in the LxF scheme plays here a key role, as is the case for the theoretical study in~\cite{dalibard2006homogenization}.

\medskip

Finally, the LxF-homogenized problem reads
\begin{equation}\label{eq: homog. problem LF 2d}
  \begin{cases}
    \partial_t \bar{u}(t,x) + \bar{b} \cdot \nabla_x \bar{u}(t,x) = 0,
    \\[4pt]
    \bar{u}(t = 0,x) = \bar{u}_0(x),
  \end{cases} 
\end{equation}
with the constant effective velocity
\begin{equation}\label{eq: effective velocity 2d}
  \bar{b} = \int_Y w \, b + \frac{c_t}{2} \, \text{div} \, (w \, b) \, b.
\end{equation}

\section{Generalities on the numerical approach} \label{sec:generalities}

We briefly explain in this section how we put in action, numerically, our strategy of approximation of the original LxF scheme. The numerical experiments themselves, along with our comments on the results, are postponed to Section~\ref{sec:numericresults}, considering dimensions $d = 1, 2$.

\medskip

Before turning to our actual tests, let us recall the exact solutions involved (these continuous solutions are denoted by lower-case letters, while their discrete approximations are denoted throughout the text with upper-case letters):
\begin{itemize}
\item $u^\veps_{c_x}$, solution to the modified equation~\eqref{eq: equiv eq-bis-2d} associated to the Lax-Friedrich scheme~\eqref{eq: lax-friedrich 2d},
\item $u^0_{c_x} = \bar{u}_{c_x} \, w_{c_x}(\cdot/\veps)$, where $\bar{u}_{c_x}$ is the solution to the equation~\eqref{eq: homog. problem LF 2d} obtained from homogenization of the modified equation~\eqref{eq: equiv eq-bis-2d}, and where $w_{c_x}$ is the corrector, solution to~\eqref{eq: cell problem 2d},
\item $\displaystyle u^0 = \lim_{\veps \to 0} u^\veps$, solution to the homogenized equation of the hyperbolic oscillatory equation~\eqref{eq: linear hyperbolic conservation law},
\end{itemize}
and $u^\veps$, solution to the oscillatory equation~\eqref{eq: linear hyperbolic conservation law}.

\medskip

The convergence diagram of Figure~\ref{fig: diagram solutions} illustrates the situation for these exact solutions.


\begin{figure}[H]
  \centering
  \begin{tikzpicture}[scale=0.9,
      >=Stealth,
      arr/.style={-{Stealth[length=2.5mm]}, thick} 
    ]


    \node (TL) at (0,0)
    {\small $\partial_t u_{c_x}^\veps + \text{div}_x \, (u_{c_x}^\veps \, B_{c_x}^\veps) = \veps \, \text{div}_x \, (A_{c_x}^\veps \nabla_x u_{c_x}^\veps) \qquad$};

    \node (TR) at (8,0)
    {\small $\qquad \partial_t \bar{u}_{c_x} + \bar{b}_{c_x} \cdot \nabla_x \bar{u}_{c_x} = 0$};

    \node (BL) at (0,-4)
    {\small $\partial_t u^\veps + \text{div}_x \, (u^\veps \, b^\veps) = 0 \qquad$};

    \node (BR) at (8,-4)
    {\small $\qquad \partial_t \bar{u}_{c_x=0} + \bar{b}_{c_x=0} \cdot \nabla_x \bar{u}_{c_x=0} = 0$};


    \draw[arr] (TL.east) -- (TR.west)
    node[midway, above] {$\veps \to 0$};

    \draw[arr] (BL.east) -- (BR.west)
    node[midway, above] {$\veps \to 0$}
    node[midway, below] {not always};

    \draw[arr, shorten >=5mm, shorten <=5mm] (TL.south) -- (BL.north)
    node[midway, right] {$c_x \to 0$};

    \draw[arr, shorten >=5mm, shorten <=5mm] (TR.south) -- (BR.north)
    node[midway, right] {$c_x \to 0$};
\end{tikzpicture}
\caption{Convergence diagram between the various continuous solutions with respect to the microscale~$\veps$ and the discretization parameter $c_x$ (for some quantities, we have changed the notation with respect to the text so as to make explicit the dependence on the parameter $c_x$). Taking $c_x \to 0$ is equivalent to taking the discretization parameters $h, \delta t \to 0$ at $\veps$ fixed.}
\label{fig: diagram solutions}
\end{figure}

\noindent To some of these exact solutions correspond the following discrete approximations:
\begin{itemize}
\item $U^\veps$, abbreviated as \emph{the LxF solution}, solution to the Lax-Friedrichs scheme~\eqref{eq: lax-friedrich 2d} applied to the oscillatory equation~\eqref{eq: linear hyperbolic conservation law} with a fine grid size, which is our reference object,
\item $U^0_{c_x}$, abbreviated as \emph{the LxF-homogenized approximation}, and approximation to~$u^0_{c_x}$, defined as $U^0_{c_x} = \bar{U}_{c_x} \, W_{c_x}(\cdot/\veps)$, where~$\bar{U}_{c_x}$ is obtained (using the Lax-Wendroff scheme) as an approximation of the solution~$\bar{u}_{c_x}$ to the homogenized equation~\eqref{eq: homog. problem LF 2d}, and where $W_{c_x}$ is a (finite element) approximation of the corrector $w_{c_x}$,
\end{itemize}
and $U^\veps_{\rm exact}$, discrete solution computed as accurately as possible (in our tests, with a Lax-Wendroff scheme and a very fine grid size), and used as ``exact'' solution for the oscillatory equation~\eqref{eq: linear hyperbolic conservation law} itself.

\medskip

We recall that the dependency upon~$c_x$ we explicitly mention above also means a dependency upon~$c_t$, since the time discretization parameter~$c_t$ is defined from~$c_x$ by~\eqref{eq: scale cfl}. In our numerical tests, we take $\alpha=0.9$ to satisfy the CFL condition.

\paragraph{Cell problem, effective velocity and preparation of the initial condition.} As a first step, we solve the cell problem~\eqref{eq: cell problem 2d} using a $\mathbb{P}_1$ Finite Element Method. The domain is the periodic cell~$Y = (0,1)^d$. It is discretized with $I = 60$ elements along each direction $1 \leq i \leq d$. In principle, the numerical solution to~\eqref{eq: cell problem 2d} could suffer from spurious oscillations because of the small diffusivity constant (of the order of the discretization parameter $c_x \ll 1$) compared to the advection term~$B$ (which is of order one). However, the fact that the boundary conditions imposed are periodic, and not Dirichlet boundary conditions, seems to eliminate the possible instabilities. We specifically checked that there is no significant improvement of our results when a stabilized (SUPG, see~\cite{brooks-hughes1982,quarteroni1994numerical}) formulation is employed. 

Once the cell problem is numerically solved, we compute the effective velocity $\bar{b}$ from~\eqref{eq: effective velocity 2d}. Likewise, we use the numerical approximation~$W_{c_x}$ of the corrector to construct a well-prepared initial condition, as the product~\eqref{eq: well-prep}, from any arbitrary macroscopic function $\bar{u}_0(x)$ (note that this computation has to be performed for each discretization parameter~$c_x$). This well-prepared initial condition is used to compute $U^\veps$, $U^\veps_{\rm exact}$, and also when we compute an as-accurate-as-possible approximation to $u^\veps_{c_x}$.

\paragraph{Solution of the LxF-homogenized problem.} For the numerical approximation $\bar{U}_{c_x}$ of the solution to the LxF-homogenized problem~\eqref{eq: homog. problem LF 2d}, we use a Lax-Wendroff finite volume scheme on a uniform spatial grid $\Omega_H$ with grid size $H$, independent of the small scale $\veps$. The time-step $\Delta T$ is given by the CFL condition for the constant effective velocity $\bar{b}$. We note that the Lax-Wendroff scheme is only used here as a representative higher-order finite volume method. Other, possibly more accurate, schemes could also be considered (see~\cite{leveque2002fv}). In any event, the equation being linear and the effective velocity being constant here (because of the periodicity assumed on the original oscillatory velocity), we could equally well explicitly solve the homogenized equation~\eqref{eq: homog. problem LF 2d} analytically. We however wish to lay some groundwork in preparation for more general settings, such as nonlinear conservation laws, or even linear cases with velocities $b(x, x/\veps)$ depending also on the slow variable, for which the same homogenized formulation applies. We then evaluate the solution $\bar{U}_{c_x}$ on the fine grid $\Omega_h$, multiply it with the solution $W_{c_x}$ to the cell problem, and finally obtain our approximation
$$
U^0_{c_x} = \bar{U}_{c_x} \, W_{c_x}(\cdot/\veps)
$$
on $\Omega_h$.

\paragraph{Metrics of comparison.} Given the multiscale nature of the problem we consider, it may be delicate to compare our numerical solutions with one another. For this purpose, we define several \emph{metrics} that we use for experiments, specifically in the one-dimensional setting (in the presentation of these metrics, we thus refer to equations in the one-dimensional setting, in contrast to the presentation above, where we referred to the two-dimensional setting). We introduce the notions of
\begin{itemize}
\item \textbf{Effective velocity}: comparing solutions to transport-like equations require to first ``renormalize'' them so that they share the same moving frame. Put differently, we need to first evaluate their effective velocities and compare these velocities, before comparing the two solutions in more details (which will be the purpose of the next two metrics). For the homogenized solution $U^0_{c_x}$, the effective velocity is $\bar{b}$ given by~\eqref{eq: effective velocity}. It parametrically depends upon the parameter $c_x$. For $\veps$ vanishing, we can assume that the effective velocity of the Lax-Friedrichs solution is also given by~\eqref{eq: effective velocity}, where $c_x$ is given by the discretization step $h$. On the other hand, an accurate approximation of the continuous solution $u^\veps$ to~\eqref{eq: linear hyperbolic conservation law} travels at the velocity $b^\star$ (given in~\eqref{eq: harmonic average b}) associated to the homogenized solution $u^0$.
\item \textbf{Relative error in the local average}: given a solution $g(t,x)$, we define its local average on each point of the grid by
  \begin{equation*}
    \hat{g}(t_n,x_i ) := \frac{1}{\veps}\int_{x_i - \veps/2 }^{x_i + \veps/2} g(t_n, s) \, ds.
  \end{equation*}
  We then compute the relative error in the local average as
  \begin{equation} \label{eq:local-average}
    \mathcal{\hat E}(g_1, g_2) := \frac{\|\hat g_1-\hat g_2\|_{L^\infty(0,T;\,L^2(\Omega))}}{\|\hat g_1\|_{L^\infty(0,T;\,L^2(\Omega))}}.
  \end{equation}
  Of course and in line with our previous item, if the two solutions compared travel at different effective velocities, they must be preliminary shifted, using~$\hat{g}(t_n,x_i + c \, t_n)$, with $c$ the effective velocity of $g(t,x)$, before the error~\eqref{eq:local-average} is evaluated.
\item \textbf{Relative error in the fine scale oscillations}. In order to eventually focus on the error in the oscillations of the solution $g(t,x)$, we subtract to any solution~$g$ its local average
  \begin{equation*} 
    \tilde{g}(t,x) := g(t, x) - \hat{g}(t,x).
  \end{equation*}
  As above, if the two solutions $g_1(t,x)$ and $g_2(t,x)$ to be compared turn out to travel at different effective velocities, we also shift each~$\tilde{g_i}$ suitably. We then measure the relative error between the oscillations of the two solutions with
  \begin{equation*}
    \mathcal{\tilde E}(g_1, g_2) := \frac{\|\tilde g_1-\tilde g_2\|_{L^\infty(0,T;\,L^2(\Omega))}}{\|\tilde g_1\|_{L^\infty(0,T;\,L^2(\Omega))}}.
  \end{equation*}
\end{itemize}
In what follows, we respectively denote $\mathcal{E}_{ave}$ and $\mathcal{E}_{osc}$ the discrete counterparts of~$\mathcal{\hat E}$ and $\mathcal{\tilde E}$.

\section{Numerical results} \label{sec:numericresults}

Our purpose in this section is to assess how accurately the LxF-homogenized solution $U^0_{c_x}$ approximates the solution $U^\veps$ to the Lax–Friedrichs scheme~\eqref{eq: lax-friedrich} (or~\eqref{eq: lax-friedrich 2d} in 2D), which is our object of reference. The error parametrically depends upon the microscale $\veps$ and the discretization parameter $c_x = h/\veps$ defining the number of grid points per microscale. More precisely, we examine the two regimes: $\veps \to 0$ with $c_x$ fixed, and $c_x \to 0$ with $\veps$ fixed, bearing in mind that in both cases $h = c_x \, \veps$ with $c_x <1$.

The first regime (in Section~\ref{subsec:numericseps1d}) illustrates the asymptotic analysis~$\veps \to 0$ of Section~\ref{subsec:homog}. In the second regime (see Section~\ref{subsec:refinementregime}), $\veps$ is held fixed and the Lax-Friedrichs discretization is refined, and we investigate how our numerical approach, based on homogenization, behaves when the scheme better and better approximates the original continuous equation. Both Sections~\ref{subsec:numericseps1d} and~\ref{subsec:refinementregime} aim at illustrating the situation in the ``phase'' space~$(\veps,c_x)$.

Schematically, the error committed when approximating the solution~$U^\veps$ of the original scheme by using the solution~$U^0_{c_x}$ of our homogenized equation may be estimated using the triangle inequality (and with self-explanatory notations)
\begin{multline} \label{eq:triangular}
  \left\|\tt{scheme} ({c_x},\veps) - \tt{homogen.modified\, eqn}({c_x}) \right\| \\ \leq \left\|\tt{scheme}({c_x},\veps)- \tt{modified\, equation}({c_x},\veps) \right\| \\ + \left\|\tt{modified\, equation}({c_x},\veps) - \tt{homogen.modified\, eqn} ({c_x}) \right\|,
\end{multline}
thereby distinguishing two sources of error. On the right-hand side of~\eqref{eq:triangular}, the first term represents how accurately the modified equation mimics the solution of the numerical scheme (a difference that essentially depends upon~$c_x$, and parametrically on~$\veps$). The second term corresponds to the homogenization error depending primarily on $\veps$, and parametrically on~$c_x$. The purpose of our numerical tests in Sections~\ref{subsec:numericseps1d} and~\ref{subsec:refinementregime} is to investigate the various regimes where either of the two terms becomes dominant.

Section~\ref{subsec:performance} finally compares quantitatively the performance of our approach with that of a classical scheme put in action on the original oscillatory equation~\eqref{eq: linear hyperbolic conservation law}, both in terms of accuracy reached at a given computational cost, and cost required to reach a given accuracy. The qualitative supremacy of the approach based on homogenization is not unexpected. Indeed, since $h = c_x \, \veps$ and $\delta t = c_t \, \veps$, we have
$$
N_x = (1/h)^d = (c_x \, \veps)^{-d} \propto \veps^{-d}, \qquad N_t = T / \delta t \propto \veps^{-1},
$$
such that the cost of computing $U^\veps$ grows as $\veps^{-(d+1)}$, whereas that of $U^0_{c_x}$, in sharp contrast, does not depend on $\veps$. 

\subsection{Asymptotic behavior for vanishing $\veps$} \label{subsec:numericseps1d}

In the following, we compare $U^0_{c_x}$ and~$U^\veps$ when decreasing the microscale $\veps$, while fixing the number of elements of the fine grid per $\veps$, thus the parameter $c_x$. In the one-dimensional setting of the present section and Section~\ref{subsubsec:refine1d}, all our simulations have been carried out over the time interval $[0,T]$ with the final time $T = 2$, given the smooth and periodic velocity field
\begin{equation}\label{eq:beps-12}
  b(x / \veps) = 2.1 - \sin(4 \pi x/\veps),
\end{equation}
and the macroscopic initial condition $\bar{u}_0 (x) = \cos(2 \pi x)$. The solution to the homogenized equation~\eqref{eq: homog. problem LF} is computed over a coarse grid $\Omega_H$ with the grid size $H = 0.005$.

\medskip

\begin{figure}[H]
  \centering
  \includegraphics[width=\textwidth]{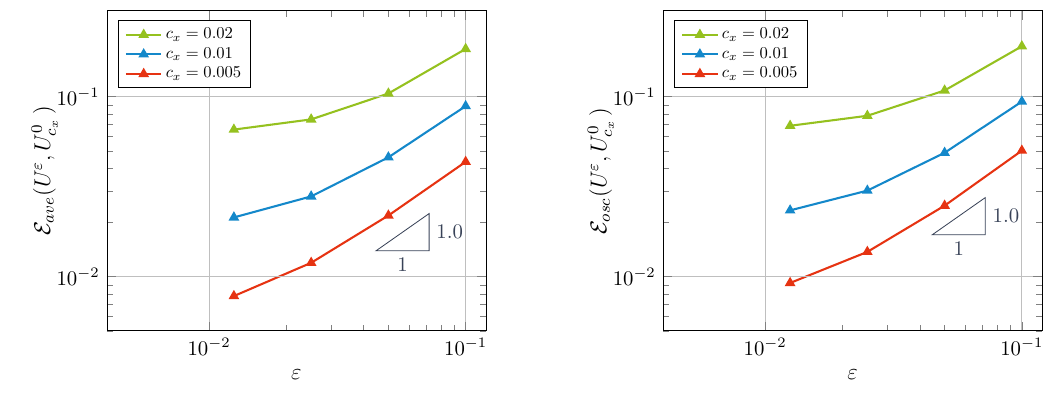}
  \caption{Error between the LxF solution and its LxF-homogenized approximation, for $\veps \in [0.01, 0.1]$, evaluated in terms of local average (left) and oscillations (right).}
  \label{fig: convergence_eps}
\end{figure}

Our results are reported in Figure~\ref{fig: convergence_eps}, showing the error between the two solutions in terms of local average and oscillations (note that we can skip the comparison in terms of our first metric, since by construction the solutions $U^0_{c_x}$ and~$U^\veps$ have an identical effective velocity), for three different $c_x \in \{0.02, \, 0.01, \, 0.005\}$.

At fixed~$c_x$, the error decreases with~$\veps$ but actually reaches a plateau for a certain, sufficiently small value of~$\veps$. This plateau is only suggested in Figure~\ref{fig: convergence_eps}, focusing on the smallest values of~$\veps$ toward the left-hand side of the graphs. Indeed, since our homogenization process has been conducted on the modified equation associated to the original Lax-Friedrichs scheme and evidently not on the scheme itself, the \textit{distance} between the modified equation and the scheme implicitly affects the error observed here. And this distance sensitively depends on the discretization parameter $c_x$, but is uniformly bounded in~$\veps$. The latter fact (which was \emph{postulated} in our derivation, following on Remark~\ref{rk:epsilon}) is specifically illustrated below in Figure~\ref{fig:lxfvsmod-two}. We note for completeness that, in this figure, the modified equation is solved using a fine $\mathbb{P}_1$ finite element discretization.

Before the plateau is reached, we observe in Figure~\ref{fig: convergence_eps} that, for such intermediate values of~$\veps$, the errors decay approximately linearly with respect to $\veps$, while for smaller values of~$\veps$ the convergence slows down.

Independently, smaller values of~$c_x$ lead to a uniform reduction of the error.

\begin{figure}[H]
  \centering
  \begin{subfigure}[t]{0.49\textwidth}
    \includegraphics[width=\textwidth]{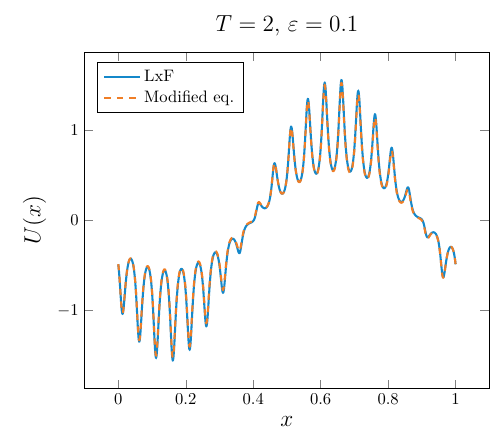}
  \end{subfigure}
  \begin{subfigure}[t]{0.49\textwidth}
    \includegraphics[width=\textwidth]{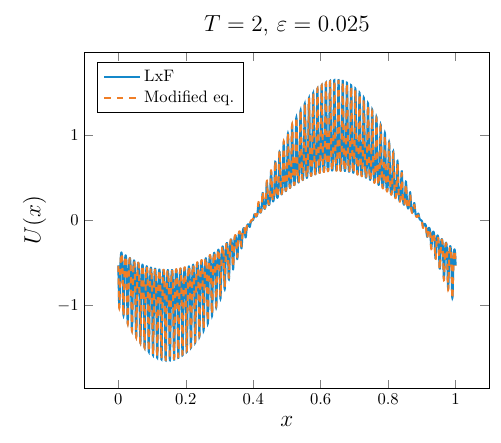}
  \end{subfigure}
  \caption{Solutions to the LxF scheme and to the modified equation at time $T = 2$, with, respectively on the left and the right hand side, the microscale~$\veps = 0.1$ and~$\veps= 0.025$ (the parameter~$c_x = 0.005$, which corresponds to 200 grid points per $\veps$, is fixed). When decreasing $\veps$, the distance between the two solutions does not vary.}
  \label{fig:lxfvsmod-two}
\end{figure}

\begin{remark}[Higher dimensional settings]
The convergence (when $\veps \to 0$) of the modified equation~\eqref{eq: equiv eq-bis-2d} towards the LxF-homogenized equation~\eqref{eq: homog. problem LF 2d} is guaranteed, in the linear setting, in any dimension and for any flow topology, following the analysis performed in~\cite{dalibard2006homogenization}. We have checked on a typical one-dimensional example of velocity field that the distance between the scheme and the modified equation remains uniformly bounded in $\veps$ (see Figure~\ref{fig:lxfvsmod-two}) and the same result is expected to hold in higher dimensions. We thus expect that, in dimensions higher than one, and similarly as in the one-dimensional setting, in the regime $\veps \to 0$ with $c_x$ fixed, the distance between the solutions to the LxF scheme and to the LxF-homogenized model decreases until reaching a lower bound, the value of which depends on the discretization parameter $c_x$.
\end{remark}

\subsection{Asymptotic behavior for $c_x$ vanishing} \label{subsec:refinementregime}

We now turn to the distance between the two solutions $U^0_{c_x}$ and~$U^\veps$, for $c_x \to 0$ (thus increasing the number of grid points per microscale), at fixed $\veps$. In this regime, the LxF solution~$U^\veps$ becomes an increasingly accurate (and expensive) approximation of the solution to the oscillatory equation~\eqref{eq: linear hyperbolic conservation law}, whereas the homogenized solution~$U^0_{c_x}$ is computed on the coarse grid $\Omega_H$ with a cost independent of $c_x$.

\subsubsection{One-dimensional setting} \label{subsubsec:refine1d}

In the one-dimensional setting, we recall that we perform our tests for the velocity field~\eqref{eq:beps-12}.

\begin{figure}[htbp]
  \centering
  \includegraphics[width=\textwidth]{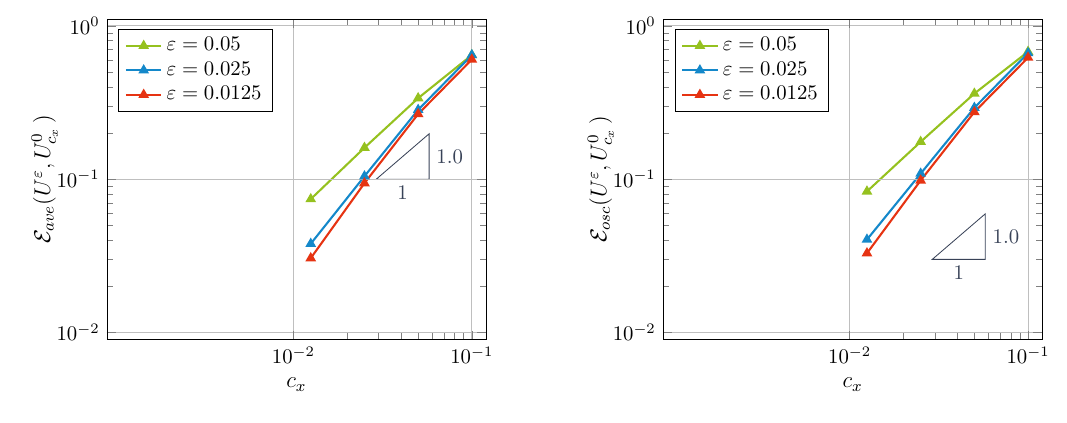}
  \caption{Error between the LxF solution and its LxF-homogenized approximation for $c_x \in [0.01, 0.1]$, evaluated in terms of local average (left) and oscillations (right).}
  \label{fig: convergence_cx}
\end{figure}

As $c_x \to 0$, the modified equation approximates the scheme more accurately. The first term in~\eqref{eq:triangular} thus decreases. We expect the error between the two solutions to correspondingly decrease. Figure~\ref{fig: convergence_cx} indeed shows the decay, both for local averages (left) and for fine-scale oscillations (right). This phenomenon is robust with respect to the value of~$\veps$.

\begin{figure}[H]
  \centering
  \includegraphics[width=0.5\textwidth]{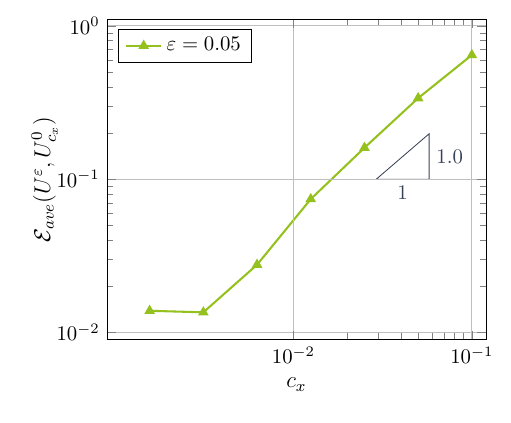}
  \caption{Error between the LxF solution and its LxF-homogenized approximation for $c_x \in [0.001, 0.1]$ and $\veps = 0.05$, evaluated in terms of local average.}
  \label{fig: convergenceeps0.05_cx}
\end{figure}

However, in the limit $c_x \to 0$, we cannot expect the error to completely vanish either, because of the second term in~\eqref{eq:triangular} that measures the homogenization error, a function of~$\veps$. This is indeed observed in Figure~\ref{fig: convergenceeps0.05_cx}, where the error reaches a lower bound close to $1\%$. 

\begin{figure}[H]
  \centering
  \begin{subfigure}[b]{\textwidth}
    \centering
    \includegraphics[width=0.49\textwidth]{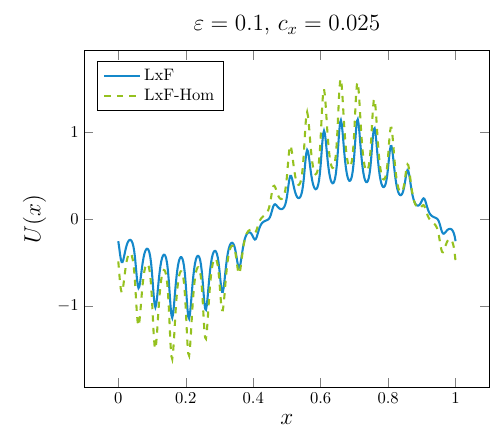}\hfill
    \includegraphics[width=0.49\textwidth]{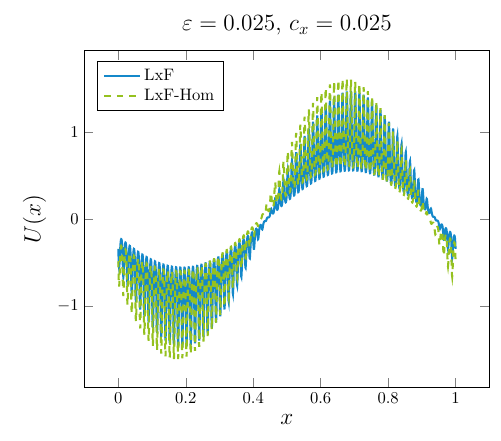}
  \end{subfigure}
  \begin{subfigure}[b]{\textwidth}
    \centering
    \includegraphics[width=0.49\textwidth]{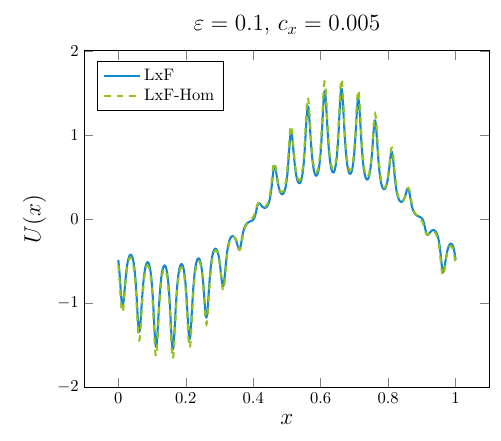}\hfill
    \includegraphics[width=0.49\textwidth]{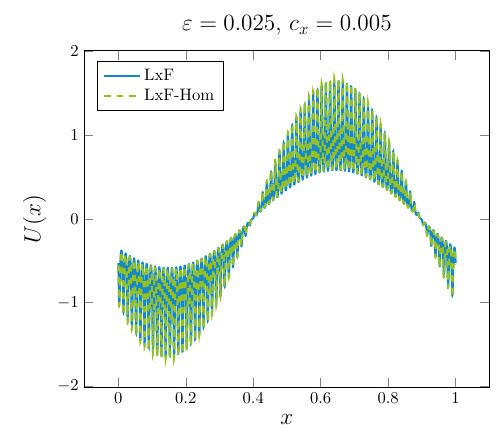}
  \end{subfigure}
  \caption{Comparison between the LxF solution and its LxF-homogenized approximation at $T = 2$, for $\veps \in \{0.1, 0.025\}$ and $c_x \in \{0.025, 0.005\}$.}
  \label{fig: 1d-simulations}
\end{figure}

In order to better visualize the respective impact of either of the two parameters $c_x$ and $\veps$ on the quality of the approximation of $U^\veps$ by $U^0_{c_x}$, we provide in Figure~\ref{fig: 1d-simulations} the profiles of these two solutions at final time, for a selection of values of~$(c_x,\veps)$.

Finally, we specifically investigate the situation for $\veps$ asymptotically small. Homogenization theory tells that, for this one-dimensional setting (and we will see later on that the question is much more delicate in higher dimensions), the original oscillatory hyperbolic equation~\eqref{eq: linear hyperbolic conservation law} homogenizes into the linear transport equation~\eqref{eq:effective-1d-lin} with constant velocity $b^\star = \langle b^{-1} \rangle^{-1}$. The numerical tests we performed for Table~\ref{tab:velocity1d} confirm that, in the limit $c_x\to 0$, the effective velocity~$\bar{b}$ appearing in~\eqref{eq: homog. problem LF} indeed converges to $b^\star$. Put differently, in this simple setting, the diagram of Figure~\ref{fig: diagram solutions} commutes.

\begin{table}[htbp]
\centering\small\setlength{\tabcolsep}{5pt}
\begin{tabular}{@{}cccccc@{}}
\toprule 
$\bar{b} \, (c_x = 0.1)$ & $\bar{b} \, (c_x = 0.05)$ & $\bar{b} \, (c_x = 2.5\text{e-2})$ & $\bar{b} \, (c_x = 1.25\text{e-2})$ & $\bar{b} \, (c_x = 6.25\text{e-3})$ & $b^\star$ \\
\midrule 
 1.91154 & 1.88232 & 1.86075 & 1.85101 & 1.84781 & 1.84662 \\
\bottomrule
\end{tabular}
\caption{Effective velocities of the LxF homogenized model, approaching the harmonic average of $b$ (i.e. the effective velocity of the hyperbolic regime) as $c_x \to 0$.}
\label{tab:velocity1d}
\end{table}

\subsubsection{Two-dimensional settings} \label{subsec:numerics2D}

We now consider settings in dimension two, with again the microscale $\veps$ fixed at a small value and a vanishing parameter $c_x$. Since the scheme is consistent, the LxF solution approaches $u^\veps$, solution to~\eqref{eq: linear hyperbolic conservation law}. We might also presume that the LxF solution is close to that of the modified equation (since $c_x \to 0$) and thus to the solution of the LxF-homogenized problem~\eqref{eq: homog. problem LF 2d} (since $\veps \ll 1$). However, the LxF homogenized problem is a pure transport problem, while the homogenized limit of~\eqref{eq: linear hyperbolic conservation law} is not necessarily of the same nature. It is therefore interesting to investigate the distance between the various solutions at play, the proximity of which may depend on the respective values of $c_x$ and $\veps$.

In the multidimensional setting, we recall that the homogenized limit of~\eqref{eq: linear hyperbolic conservation law} is not explicitly known, in general. When an effective model exists, its structure depends strongly on the topology of the flow. It may be a transport equation, an equation involving non-local effects, or even a system of equations. In contrast, the modified equation~\eqref{eq: equiv eq-bis-2d} associated to the Lax-Friedrichs scheme contains numerical viscosity. The effective equation is a transport equation, \eqref{eq: homog. problem LF 2d}, independently of the flow topology.

\medskip



In this section, we test our approach on a selection of velocity fields (see Table~\ref{tab:velocity2d}) exhibiting different flow structures. For the first three test cases, our approach turns out to be an effective approximation of the LxF solution~\eqref{eq: lax-friedrich 2d} in the limit $c_x \to 0$. The shear-flow case (Test Case 4 in Table~\ref{tab:velocity2d}) is of different nature and highlights the limitations of our approach. 

\begin{table}[htbp!]
\centering\small\setlength{\tabcolsep}{5pt}
\begin{tabular}{@{}clccl@{}}
\toprule
 & $b(y_1,y_2)$ & $\text{div} \, b = 0$ & theory \\
\midrule
1 & $[\sin(2\pi y_2) + \sqrt{2}, \ \cos(2\pi y_1) + \pi/2]$ & \cmark & \cmark \\
2 & $[\sin(2\pi y_2) + 1.1, \ 2.1 - \cos(2\pi y_1)]$ & \cmark & \textbf{?} \\
3 & $[0.1 \, \sin(2\pi y_1) \, \cos(2\pi y_2) + 0.5, \ 2.1 - \sin(2\pi y_1) \, \cos(2\pi y_2)]$ & \xmark & \textbf{?} \\
4 & $[\sin(2\pi y_2) + 1.1, \ 0]$ & \cmark & \cmark \\
\bottomrule
\end{tabular}
\caption{The velocity fields we consider are given by $b^\veps(x) = b(y)$ with $y=x/\veps$, for the functions $b$ given here.}
\label{tab:velocity2d}
\end{table}

\medskip

\begin{figure}[H]
  \centering
  \includegraphics[width = 0.4\textwidth]{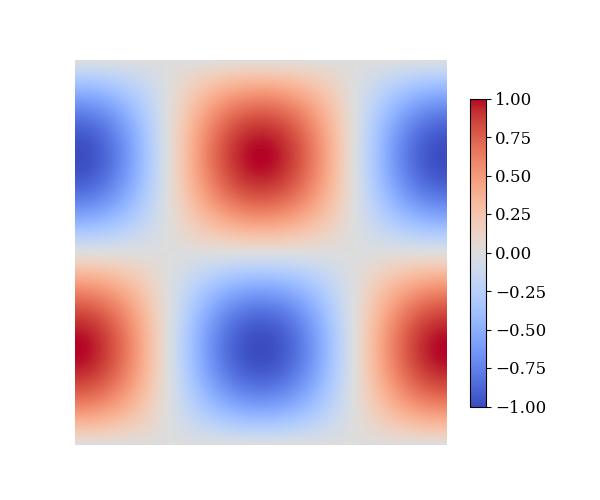}
  \caption{Initial condition~\eqref{eq:IC_2d}.}
  \label{fig:IC_2d}
\end{figure}

\begin{figure}[htbp!]
  \centering
  \includegraphics[trim={1cm 0 1cm 0},width=0.49\textwidth]{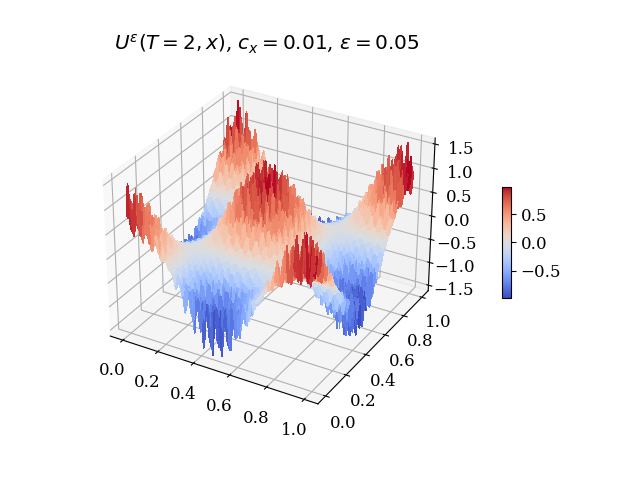} \hfill
  \includegraphics[trim={1cm 0 1cm 0},width=0.49\textwidth]{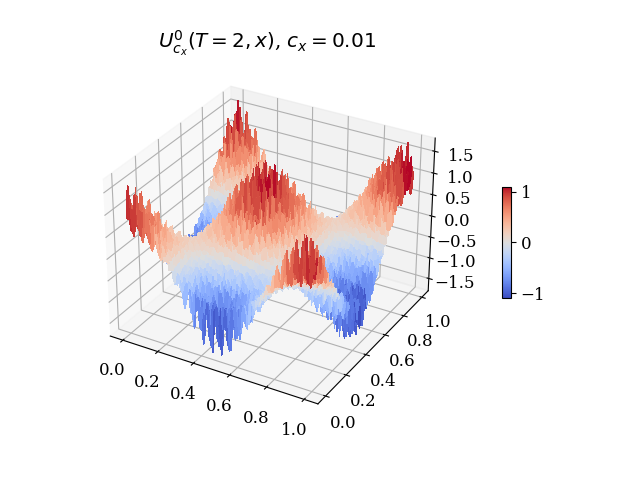}
  \caption{Comparison between the LxF solution (left) and its LxF-homogenized approximation (right) at $T = 2$ ($\veps = 0.05$, $c_x = 0.01$, Test Case 3 of Table~\ref{tab:velocity2d}).}
  \label{fig: 2d-simulations}
\end{figure}

\begin{figure}[hbtp!]
  \centering
  \begin{subfigure}[t]{0.49\textwidth}
    \includegraphics[trim={1cm 0 1cm 0},width=\textwidth]{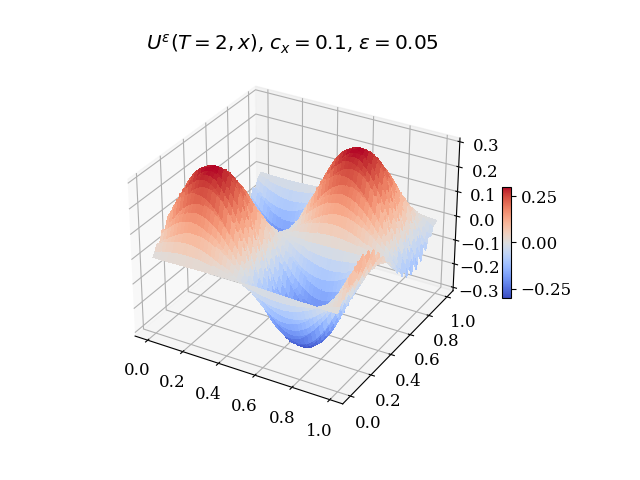}
    \caption{}
    \label{subfig: shear-simul_a}
  \end{subfigure}
  \hfill
  \begin{subfigure}[t]{0.49\textwidth}
    \includegraphics[trim={1cm 0 1cm 0},width=\textwidth]{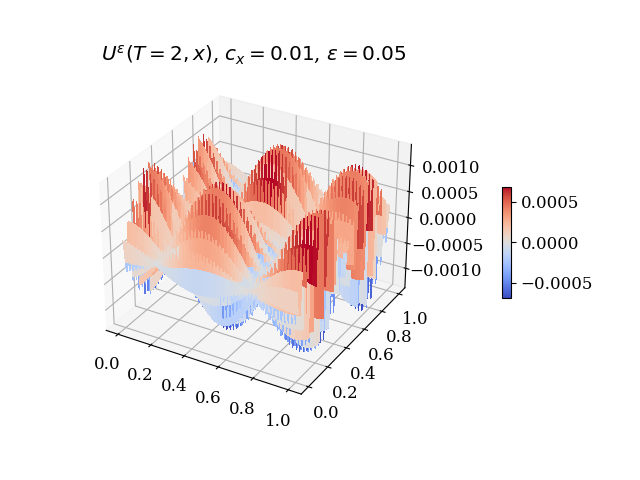}
    \caption{}
    \label{subfig: shear-simul_b}
  \end{subfigure}
  \begin{subfigure}[t]{\textwidth}
    \centering
    \includegraphics[trim={1cm 0 1cm 0},width=0.54\textwidth]{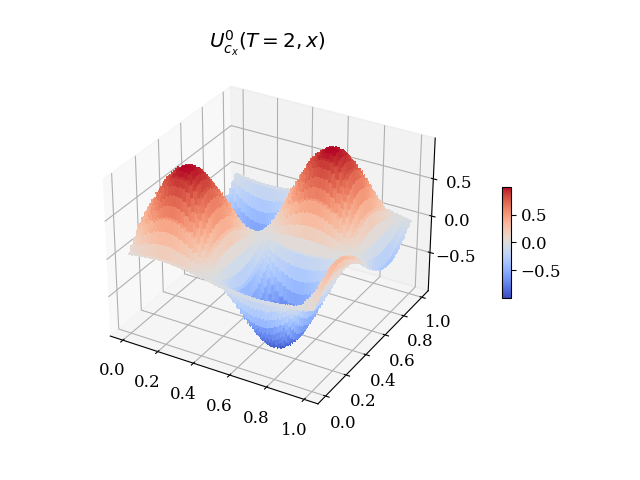}
    \caption{}
    \label{fig:shear-LxF-hom}
  \end{subfigure} 
  \caption{Comparison between the LxF solution (Figures~(a) and~(b), with two different grid sizes) and its LxF-homogenized approximation (Figure~(c)), which is independent of $c_x$ ($T = 2$, $\veps = 0.05$, Test Case 4 (shear flow) of Table~\ref{tab:velocity2d}).}
  \label{fig:shear}
\end{figure}

The following simulations have been performed for $\veps = 0.05$, over the time interval $[0,T]$ with the final time $T = 2$. The macroscopic initial condition, shown in Figure~\ref{fig:IC_2d} (and to be compared with the solutions at time $T=2$ shown in Figures~\ref{fig: 2d-simulations} and~\ref{fig:shear}), is
\begin{equation} \label{eq:IC_2d}
  \bar{u}_0 (x_1, x_2) = \cos(2\pi x_1) \, \sin(2\pi x_2).
\end{equation}

In the present two-dimensional setting, we only evaluate the relative error
\begin{equation} \label{eq:rel_error}
\frac{\| U^\veps - U^0_{c_x}\|_{L^\infty(0,T,L^2(\Omega))}}{\| U^\veps \|_{L^\infty(0,T,L^2(\Omega))}}
\end{equation}
between the LxF solution $U^\veps$ and its LxF-homogenized approximation $U^0_{c_x}$ and do not differentiate between the contributions of the local average and the oscillations. The LxF solution has been computed with a grid size $h = c_x \, \veps \in [0.0005, 0.005]$ (we recall that $\veps = 0.05$ is fixed), while the LxF-homogenized problem~\eqref{eq: homog. problem LF 2d}, the parameters of which depend on $c_x$, is discretized with the fixed step $H=0.025$.

\paragraph{Test Cases 1 to 3.} We successively consider two incompressible ($\text{div} \, b = 0$) and one compressible velocity fields. Our first test case is an incompressible flow that is ergodic (see~\cite{weinan1992homogenization,hou1992homogenization} for details on this theoretical property). It can then be shown that the function $w \equiv 1$ is the unique periodic solution to the problem $\text{div} \, (w \, b) = 0$ with $\langle w \rangle = 1$. From this property, it follows that the homogenized limit of~\eqref{eq: linear hyperbolic conservation law} is a \textit{transport equation} of the form~\eqref{eq:effective-1d-lin} with an effective velocity given as the arithmetic average of $b$ over the periodic cell.

For this incompressible velocity field, the cell problem~\eqref{eq: cell problem 2d} associated with our (\textit{parabolic}) modified equation is elliptic and has for unique solution the corrector $w \equiv 1$. As a consequence, the LxF-homogenized problem~\eqref{eq: homog. problem LF 2d} is a transport equation with an effective velocity $\bar{b}$ given as $\langle b \rangle$, independently of the value of $c_x$. We hence see that the homogenized limit of~\eqref{eq: linear hyperbolic conservation law} coincides with the LxF-homogenized problem: stated otherwise, the diagram of Figure~\ref{fig: diagram solutions} commutes. 

\begin{figure}[htbp!]
  \centering
  \includegraphics[width = 0.5\textwidth]{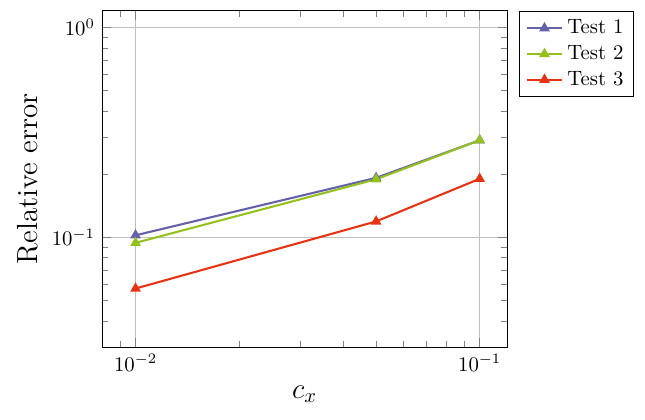}
  \caption{Relative error~\eqref{eq:rel_error} between the LxF solution and its LxF-homogenized approximation for $c_x \in [0.01, 0.1]$, with $\veps = 0.05$ and $T = 2$.}
  \label{fig: 2d-cxto0-converg}
\end{figure}

In this Test Case 1, numerical results show that, as expected, the distance between the LxF solution and its LxF-homogenized approximation decreases when $c_x$ decreases, as shown in Figure~\ref{fig: 2d-cxto0-converg} (we however expect a plateau to be reached when $c_x$ vanishes, corresponding to the homogenization error, a function of $\veps$).

\medskip


%

In the case of Test Cases 2 and 3, the homogenized limit of the hyperbolic oscillatory equation~\eqref{eq: linear hyperbolic conservation law} is, to the best of our knowledge, unknown theoretically. However, the numerical results we obtain are consistent with the fact that this limit is, like in Test Case 1, a transport equation with constant velocity. In other words, for these Test Cases 2 and 3, the diagram of Figure~\ref{fig: diagram solutions} also commutes. 

The numerical results obtained for Test Case 2 and Test Case 3 are similar to those of Test Case 1 (see Figure~\ref{fig: 2d-cxto0-converg}). The error decreases below $10\%$ for $c_x \leq 0.01$. In order to provide some visual results, we compare in Figure~\ref{fig: 2d-simulations} (for the latter case, i.e. Test Case 3) the solutions at the final time $T=2$. The LxF solution seems there slightly more diffusive than its LxF-homogenized approximation, even though a small grid step ($h = 5\text{e-4}$) has been used (the maximal value of $U^0_{c_x}$ is 1.5, while it is a bit smaller for $U^\veps$).

\paragraph{Test Case 4.} For the first three test cases, our approach not only provides an approximation to the LxF solution~\eqref{eq: lax-friedrich 2d}, but also, in the limit $c_x \to 0$, an approximation of the solution to the original oscillatory equation. This is expected from theory for Test Case 1, which is an ergodic flow. For Test Cases 2 and 3, we have no theoretical guarantee that the diagram commutes, but we do observe, based on our numerical tests, that it does.

Our fourth and final test case, Test Case 4, is a shear flow. Specifically (see the expression of $b$ in Table~\ref{tab:velocity2d}), transport occurs in the direction $x_1$ while the velocity oscillates in the transverse direction $x_2$. For this test case, we happen to know from theory (see~\cite{tartar1989,mascarenhas1993memory,AMIRAT1989,amirat1991convdiff}) that the homogenized limit of~\eqref{eq: linear hyperbolic conservation law} is a transport equation with an additional term, namely a non-local diffusion with memory. Therefore, it cannot be the case that our approach is correct in the limit $c_x \to 0$. The present final test aims at illustrating the limitation of our approach. More precisely, our approach has been derived from an homogenization process ($\veps \to 0$) at $c_x$ fixed. In the plane of parameters $(c_x,\veps)$, different asymptotic regimes can be considered. When the diagram of Figure~\ref{fig: diagram solutions} commutes, no particular caveat is expected, since the same limit is obtained, however $c_x$ and $\veps$ vanish. This is the situation of our first three test cases. Difficulties regarding the respective role of $c_x$ and $\veps$ may only be observed when the diagram does not commute. This is the situation of the present fourth test case, which we have chosen to be the shear flow, because for this precise case, we know from theory what the limits are.

\begin{figure}[H]
  \centering
  \includegraphics[width=0.5\textwidth]{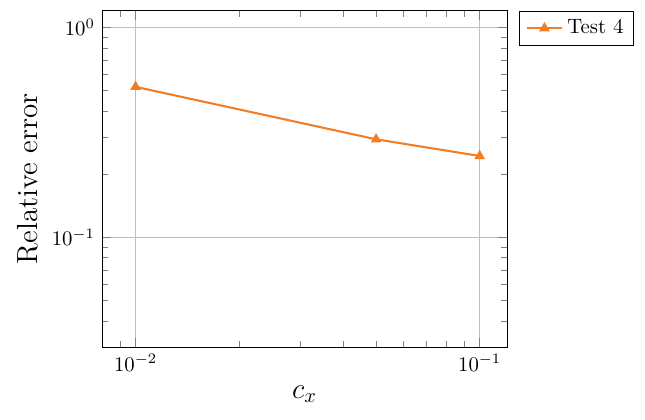}
  \caption{Relative error~\eqref{eq:rel_error} between the LxF solution and its LxF-homogenized approximation for $c_x \in [0.01, 0.1]$, with $\veps = 0.05$ and $T = 2$, in the case of shear flow (Test Case 4 of Table~\ref{tab:velocity2d}).}
  \label{fig:shear-cxto0}
\end{figure}

\begin{figure}[htbp!]
  \centering
  \includegraphics[width=0.5\textwidth]{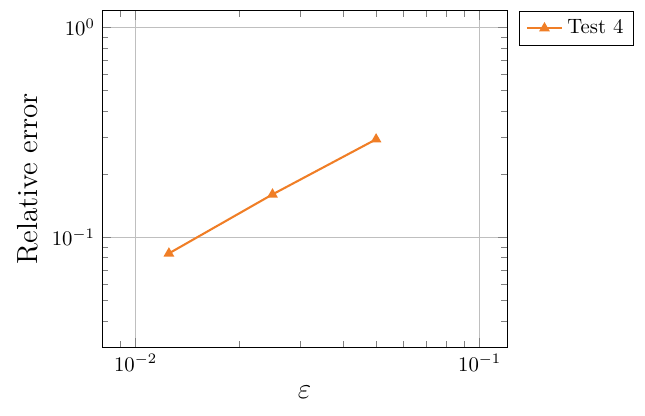}
  \caption{Relative error~\eqref{eq:rel_error} between the LxF solution and its LxF-homogenized approximation for $\veps \in [0.0125, 0.05]$, with $c_x = 0.05$ and $T = 2$, in the case of shear flow (Test Case 4 of Table~\ref{tab:velocity2d}).}
  \label{fig:shear-epsto0}
\end{figure}

Similarly as for the other tests above, we compute the numerical solution of the LxF scheme~\eqref{eq: lax-friedrich 2d} for a fixed microscale $\veps = 0.05$, letting $c_x \to 0$. The velocity field being incompressible, the unique solution to the corrector equation~\eqref{eq: cell problem 2d} is again $w \equiv 1$ and the LxF-homogenized problem therefore does not depend on $c_x$: it is a transport equation with the effective velocity $\bar{b} = \langle b \rangle$. Changing the value of $c_x$ thus only affects the LxF solution. In contrast to Test Cases 1 to 3 considered above, we observe here that the error between the LxF solution and its LxF-homogenized approximation \textit{increases} as the grid is refined, that is $c_x \to 0$ (see Figure~\ref{fig:shear-cxto0}). This of course does not contradict our analytical derivations, which have performed in the regime when $\veps \to 0$ with $c_x$ fixed. This suggests, on the other hand, that the error we commit critically depends on $c_x$, for instance as $\veps/c_x$. We may reinstate the convergence upon letting $\veps$ go to 0 at $c_x$ fixed (see Figure~\ref{fig:shear-epsto0}).

Figure~\ref{fig:shear} helps visualizing on the solutions themselves the phenomenon of Figure~\ref{fig:shear-cxto0}: on the coarsest grid ($c_x = 0.1$), the LxF solution is close to the LxF-homogenized one (compare Figure~\ref{subfig: shear-simul_a} with Figure~\ref{fig:shear-LxF-hom}), and departs from it as $c_x$ decreases (compare Figure~\ref{subfig: shear-simul_b} with Figure~\ref{fig:shear-LxF-hom}). 

%

\begin{remark} \label{rem:scaling}
In our homogenization approach, the scaling regime of the space discretization parameter $h$ in function of the microscale $\veps$ affects the prefactor of the viscous term of the modified equation~\eqref{eq: equiv eq-bis-2d} and consequently determines the homogenized problem~\eqref{eq: homog. problem LF 2d}. The specific choice we made, namely $h = c_x \, \veps$, allowed us to benefit from the theoretical findings from~\cite{dalibard2006homogenization}, which then support our formal derivation of the homogenized equation. A different scaling of the discretization parameters with respect to the microscale $\veps$ could yield a different limit problem. In particular, in the case of incompressible flows, taking $h = c_x \, \veps^2$ gives a parabolic equation of the form
$$
\partial_t u^\veps + \text{div} \, (u^\veps \, b^\veps) = \veps^2 \, \text{div} \, (A^\veps \nabla u^\veps)
$$
as modified equation (with $A^\veps = A(\cdot/\veps)$ and $A$ given by~\eqref{eq:AB-2d}), intuitively closer to the original hyperbolic equation~\eqref{eq: linear hyperbolic conservation law} than our modified equation~\eqref{eq: equiv eq-bis-2d} (and thus more likely to better approximate the purely hyperbolic regime). In the case of a shear flow (such as the one considered above) of the form $b^\veps(x_1,x_2) = b_1(x_2/\veps)$ with $b_1$ bounded away from 0, it is possible to (at least formally) derive, using a two-scale expansion, the homogenized problem that reads as a non-local diffusion equation with a memory term. This is not unexpected since it is known theoretically (since the works of Tartar and other authors, see \emph{e.g.}~\cite{tartar1989}) that, in the pure hyperbolic regime, the shear flow problem actually homogenizes into an equation with such a memory term. Investigating different scalings for the law of the space grid size $h$ in function of $\veps$ could be a possible track for further research and derivations of homogenized equations and schemes in the spirit of the present study (see~\cite{giulia_in_prep}). A related observation was made in the context of non-local conservation laws for traffic models in~\cite{colombo2021}. Numerical results obtained with the LxF scheme (in the limit when $\veps \to 0$) were more in accordance with analytical results (non-local limit) when the space grid size $h$ was of the order of $\veps^2$ than those obtained with $h$ of the order of $\veps$.
\end{remark}

\subsection{Performance assessment} \label{subsec:performance}


We illustrate in this final section the performance of our approach on a \textit{one-dimensional} test case, comparing it to a finite difference scheme put in action on the original oscillatory equation~\eqref{eq: linear hyperbolic conservation law}. For our test, we consider the velocity field $b(x/\veps) = \sin(2 \pi x/\veps) + 1.1$ and take~$\veps = 0.01$. We compare (i) the accuracy obtained for the same computational cost and (ii) the computational cost needed to reach a certain accuracy. Our \textit{reference} here is the solution $U^\veps_{\rm exact}$ to the oscillatory equation~\eqref{eq: linear hyperbolic conservation law} computed with a Lax-Wendroff scheme highly refined ($\hat{h} = \veps/100$, and correspondingly $\widehat{\delta t}$ from~\eqref{eq:cfl-general} with $\alpha=0.9$), and thus seen as an ``infinitely'' accurate approximation to $u^\veps$.

\medskip
 
Let us remark that any numerical scheme applied to the homogenized equation~\eqref{eq: homog. problem LF} is inexpensive, and even free of cost, but possibly deviates from the scheme~\eqref{eq: lax-friedrich}, unless $c_x$ is chosen sufficiently small (see the results of Section~\ref{subsec:refinementregime}). On the other hand, any scheme directly applied on the original oscillatory equation~\eqref{eq: linear hyperbolic conservation law} can be as accurate as possible, simply by decreasing the grid size $h$, but with a corresponding increase of the computational cost. Within the homogenized approach (respectively, the standard approach), we hence manipulate the model parameter $c_x$ (respectively, the discretization step $h$) to adjust for accuracy. In the specific case when we want to assign a computational cost to the two approaches, and in particular to the former one, we are led to \textit{numerically} solve a transport equation with constant velocity, which is unnecessary, since the equation has an analytical solution. Still, this procedure allows to put ourselves in a more general framework where, for example, the effective velocity $\bar{b}$ depends on the slow variable $x$ (as opposed to being constant), such that the homogenized equation indeed needs a scheme to be solved. Therefore, the goal of this final section is to go beyond the qualitative conclusions just mentioned, and to determine \textit{quantitatively} how one approach is better than the other, in terms of computational cost or accuracy.

\medskip
 
For our comparison, we use the \textit{Lax-Friedrichs scheme}~\eqref{eq: lax-friedrich}, not anymore as the reference, but as the \textit{prototype} of a classical finite difference scheme. Although some more accurate schemes could be employed, the behavior discussed below is expected to be essentially independent of this choice. In particular, any classical finite difference scheme must resolve the microscale $\veps$, requiring $h \ll \veps $, which leads to a significant computational cost. In contrast, the approach we propose is designed to operate in a coarser discretization regime, with grid sizes $h > \veps$. 

\medskip

In order to compare the methods \textit{at fixed computational cost}, we compute the solution to the LxF-homogenized problem~\eqref{eq: homog. problem LF} and the Lax-Friedrichs solution~\eqref{eq: lax-friedrich} with the \textit{same discretization step}~$h$. Unlike in Sections~\ref{subsec:numericseps1d} and~\ref{subsec:refinementregime}, the parameter $c_x$ governing the cell problem~\eqref{eq: cell problem w1} and the effective velocity~\eqref{eq: effective velocity} is here decoupled from the grid size $h$ and acts only as a regularization parameter, so that it may be chosen far smaller than the microscale would allow in a direct computation (in practice, we fix $\overline{c_x} = 1\text{e-4}$). We then evaluate the error of the LxF-homogenized solution and of the Lax-Friedrichs solution with respect to the reference $U^\veps_{\rm exact}$, in terms of local average and oscillations.

\begin{figure}[H]
  \centering
  \includegraphics[width=\textwidth]{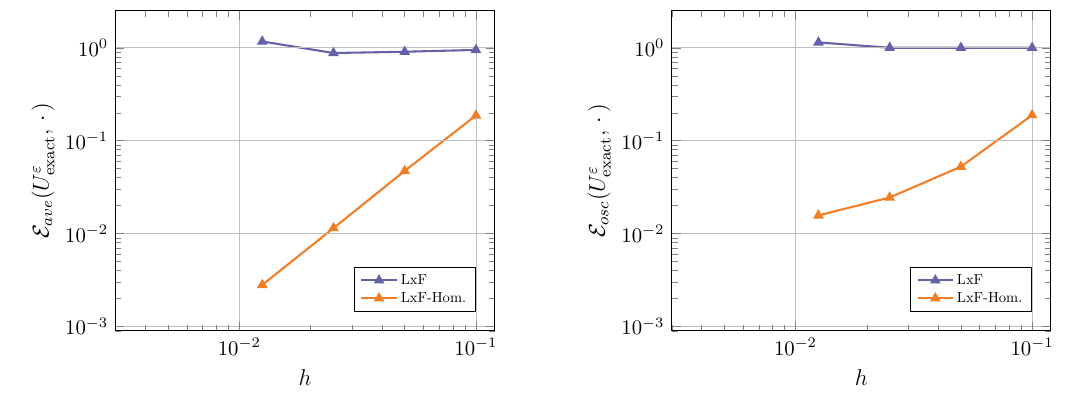}
  \caption{Accuracy of the LxF solution (blue) and of the LxF-homogenized solution (orange) at fixed computational cost: for the same discretization step $h$, we evaluate the relative error in terms of local average (left) and oscillations (right). The parameter of the homogenized problem is fixed at $\overline{c_x} = 1\text{e-4}$.}
  \label{fig: performance CPU}
\end{figure}

Numerical results, given in Figure~\ref{fig: performance CPU}, show that, for a discretization step $h \in [0.01, 0.1]$, which is larger than the microscale $\veps = 0.01$, the LxF scheme fails to capture the oscillatory solution, with relative errors close to 100\% in both the averaged and oscillatory components. In sharp contrast, the LxF-homogenized solution has a local average error (resp. error in terms of oscillations) decreasing at second order (resp. first order) in terms of $h$, and varying between 10\% and 0.1\% (resp. between 10\% and 1\%).

\begin{figure}[H]
  \centering
  \includegraphics[width=\textwidth]{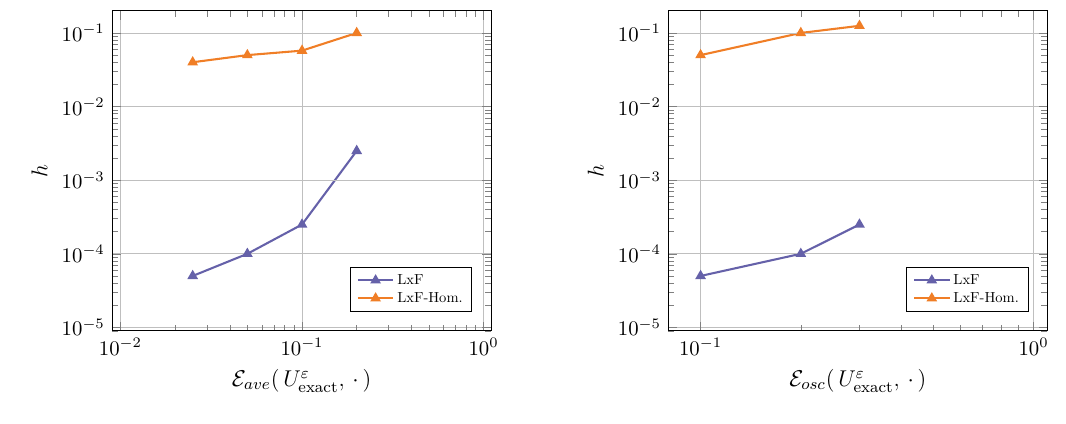}
  \caption{Computational cost of the LxF solution (blue) and the LxF-homogenized solution (orange) at fixed accuracy: evaluation of the minimum space step $h$ to obtain the same relative error in terms of local average (namely $\mathcal{E}_{ave} = 0.025$, 0.05, 0.1 or 0.2; left) or of oscillations (namely $\mathcal{E}_{osc} = 0.1$, 0.2 or 0.3; right). The parameter of the homogenized problem is fixed at $\overline{c_x} = 1\text{e-4}$.}
  \label{fig: performance accuracy}
\end{figure}

We next compare in Figure~\ref{fig: performance accuracy} the \textit{computational cost} needed by either numerical approach to obtain \textit{the same accuracy} in approximating the reference solution, both in terms of local average $\mathcal{E}_{ave} \in \{ 0.025 , \, 0.05 , \, 0.1, \, 0.2 \}$ and oscillations $\mathcal{E}_{osc} \in \{ 0.1, \, 0.2, \, 0.3 \}$. For both errors, we observe that, at fixed accuracy, the LxF-homogenized approximation captures the oscillatory solution with a space step $h$ up to three orders of magnitude larger than the discretization step needed by the Lax-Friedrichs scheme.

\paragraph{Acknowledgments.} The first two authors are grateful to Rémi Abgrall for stimulating discussions. The first author also acknowledges interesting discussions with Antoine Llor on numerical schemes for fluid dynamics. The work is supported by the EOARD under Grant Contract FA 8655-24-1-7057. The research of the first two authors is also partially supported by the ONR under Grant Contract N00014-25-1-2299. In its early stages, this work has also greatly benefited from the contributions of Amandine Boucart during her post-doc, funded by a previous ONR grant.

\bibliographystyle{plain}


\bibliography{biblio_llm.bib}

\end{document}